\documentclass[11pt]{article}

\usepackage{amsfonts,graphicx}

\newcommand{\ex}{{\bf\sf E}}               

\newcommand{\bn}{{n}}               

\newcommand{\call}{{\cal L}}

\newcommand{\calk}{{\cal K}}
\newcommand{\calr}{{\cal R}}

\newcommand{\g}{\lambda}               
\newcommand{\G}{\Lambda}               

\newcommand{\imp}{\Rightarrow}           

\newcommand{\startb}{\parindent0pt\bf}  

\newtheorem{thm}{Theorem}
\newtheorem{lem}[thm]{Lemma}
\newtheorem{pro}[thm]{Proposition}
\newtheorem{cor}[thm]{Corollary}

\def\nd{\quad {\rm and}\quad}

\def\ind{{\bf 1}}

\def\qed{\square}

\begin{document}

\title{{\LARGE\bf Heavy-Traffic Optimality of
a Stochastic Network\\
under Utility-Maximizing Resource Control}}

\author{
Heng-Qing Ye\thanks{Supported in part by the grant R-314-000-061-112
of National University Singapore.}
\\
Dept of Decision Science, School of Business\\
National University of Singapore, Singapore\\
\\
David D.\ Yao\thanks{Supported in part by NSF grant CNS-03-25495, and
by HK/RGC Grant CUHK4179/05E.
}\\
Dept of Industrial Engineering and Operations Research\\
Columbia University, New York, USA\\
}
\date{January 2006}

\maketitle

\begin{abstract}
We study a stochastic network that consists of a set of servers
processing multiple classes of jobs. Each class of jobs requires a
concurrent occupancy of several servers while being processed, and
each server is shared among the job classes in a head-of-the-line
processor-sharing mechanism. The allocation of the service
capacities is a real-time control mechanism: in each network state,
the control is the solution to an optimization problem that
maximizes a general utility function. Whereas this resource control
optimizes in a ``greedy'' fashion, with respect to each state, we
establish its asymptotic optimality in terms of (a) deriving the
fluid and diffusion limits of the network under this control, and
(b) identifying a cost function that is minimized in the diffusion
limit, along with a characterization of the so-called fixed point
state of the network.

\end{abstract}

\newpage

\baselineskip 15pt

\section{Introduction}

We study  a class of stochastic networks with concurrent occupancy
of resources, which, in turn, are shared among jobs. For example,
streaming a video on the Internet requires bandwidth from all the
links that connect the source of the video to its destination; and
the capacity (bandwidth) of each link is shared,
 according to a
pre-specified protocol, among all source-destination connections
that involve this link (including this video streaming).
A second example is a multi-leg flight on an
airline reservation system: to book the flight, seats on all legs
must be committed simultaneously.
Yet another example is a make-to-order (or, assemble-to-order)
system: the arrival of an order will trigger the simultaneous
production of all the components that are required to
configure the order.

The engineering design of Internet protocols is a center piece of
the bandwidth-sharing networks, the first example above. A key issue
here is the real-time allocation of each link's capacity to each job
class, which takes the form of solving an optimization problem for
each network state. This real-time allocation scheme makes it very
difficult if not intractable to evaluate the system performance over
any extended period as opposed to real time. For instance, if the
network is modeled as a continuous-time Markov chain, then the
transition rate from each state is itself a solution to an
optimization problem. Consequently, if one wants to optimize the system
performance over a relatively long planning horizon, it is not clear
what the resource allocation scheme one should use. In other words,
there is a gap between real-time resource control and
performance optimization in such networks.

The objective of this paper is two fold. First,
we want to overcome the intractability of
performance evaluation under dynamic resource control by developing
fluid and diffusion models: through suitable scaling of time and space, we
show that the processes of interest in the network, such as queue
lengths and workloads, converge to deterministic functions
or reflected Brownian motions, under a
broad class of utility-maximizing resource control schemes.

Our other objective is to establish the connection between the
real-time resource allocation scheme and its impact on the system
performance. We show that while the utility-maximizing control
optimizes in a ``greedy'' fashion, with respect to each state of the
network, it will, in fact, minimize a certain cost function of
performance measures such as queue lengths and workloads in the
limiting regimes of the network. In this sense, the
utility-maximizing control is asymptotically optimal.

The fluid and diffusion limits are useful in characterizing the
network performance for a very broad and important class of
resource-control policies, whereas other analytical approaches
appear to be intractable. Furthermore, the precise relationship
between the utility function (which the resource control maximizes
in each network state) and the cost function (which is minimized in
the limiting regime) provides a useful means in deciding what type
of controls (or protocols) to use so as to achieve certain
cost-related performance objectives. Another advantage of these
results is their generality, for instance, in requiring only very
mild conditions on arrival and service mechanisms.

The remainder of the paper is organized as follows. We start with an
overview of the related literature in the rest of this introductory
section. In Section \ref{sec-model}, we present details of the
network model and the utility-maximizing allocation scheme. We also
bring out a cost minimization problem, and demonstrate its close
connection with the utility maximization problem. In Section
\ref{sec-fluid}, we study the limiting behavior, in terms of fluid
scaling, of the network under the utility-maximizing resource
control, and establish the important property of {\it uniform
attraction}: the fluid limit will converge uniformly over time to
the so-called fixed-point state, which is a solution to the cost
minimization problem.
In Section \ref{sec-diffusion}, we consider diffusion scaling and
prove the main result, Theorem \ref{thm-opt-diff}: the diffusion
limit of the workload is a reflected Brownian motion, and the state
(queue-length) process converges to a fixed point that is a
minimizer of the cost objective; furthermore, both the workload and
the cost objective are minimized under the utility maximizing
control. Theorem \ref{thm-opt-diff} requires a key condition, that
there is a single bottleneck link/server in the network. (A
bottleneck link is one whose capacity is saturated by the
 offered traffic load.)
We justify this condition by showing that it is equivalent
to the so-called resource pooling condition, which is required
in many models of heavy-traffic limit.

\subsubsection*{Literature Review}

Circuit-switched technology, or, long-distance telephony, predates
the Internet protocol. The corresponding stochastic network model
 is often referred to
as the ``loss network'' --- calls that arrive when all circuits
are occupied will be blocked and lost.
Under Markovian assumptions,
the loss network in steady state
 has a product-form distribution, which, however, is still
computationally intractable due to the
combinatorial explosion of the state space;
refer to Ross \cite{keith}.

Whitt's pioneering work \cite{whitt} called
attention to concurrent resource occupancy. It was motivated
by studying the blocking phenomenon in loss networks.
Kelly \cite{kellylarge,kellyloss} developed a fixed-point approximation
for the blocking probabilities in loss networks (i.e.,
the blocking probabilities are solutions to a nonlinear system
of equations), and studied the
properties of the fixed-point mapping.

As technology evolves, research in loss networks has moved into a
new area, so-called bandwidth-sharing networks, where the service
capacity (bandwidth) on each link is shared at any time among all
the jobs in process at the link. These networks  not only model
existing protocols on the Internet (e.g., TCP), but also lead to
studies of new bandwidth allocation schemes. Most of these real-time
allocation schemes try to enforce some notion of fairness. Examples
include the max-min fair allocation of Bertsekas and Gallager
\cite{bg}, and its extensions to proportional fairness and other
related notions in Fayolle {\it et al.}~\cite{fayolle}, Kelly
\cite{kelly1,kelly2}, Kelly {\it et al.}~\cite{kellyetal}, Low
\cite{Low02}, Massoulie and Roberts \cite{mr}, Mo and Walrand
\cite{mw}, Wen and Arcak \cite{WenArcak04}, and Ye \cite{y1} among
many others. Typically, the real-time capacity allocation takes the
form of a solution to an optimization problem, with the objective
being a utility function (of the state and the allocation), and the
constraints enforcing the capacity limit on the links. Even under
Markovian assumptions,
this results in a continuous-time Markov chain with state-dependent
transition rates that are the allocation schemes resulted from
solving the above optimization problem. There is no analytically
tractable way to capture the probabilistic behavior of such a Markov
chain, such as its steady-state distribution and related performance
measures.

This has motivated research on fluid models of such networks, where
the discreteness of jobs and randomness of their arrival and service
mechanisms are transformed, via law-of-large-numbers scaling, into
continuous flows traversing the network in a purely deterministic
manner. Refer to for example Bonald and Massoulie \cite{bm}, de
Veciana {\it et al.}~\cite{vlk}, Kelly and Williams \cite{kw},
Massoulie and Roberts \cite{mr}, and Ye {\it et al.}~\cite{ye}.
Fluid models turn out to be particularly effective in studying the
stability of bandwidth-sharing networks under various capacity
allocation schemes.
 (This is notably an extension of
earlier studies on the stability of traditional queueing networks
under various service disciplines. Refer to, e.g., Dai
\cite{dai}.)
 For instance, it is shown in \cite{ye} that under the usual traffic
condition (i.e., each link has enough capacity to process all the
traffic that uses the link), many ``fair'' allocations and related
utility-maximizing schemes result in a stable network, whereas
certain priority and (instantaneous) throughput-maximizing schemes
may lead to instability.

The fluid model in \cite{kw} provides an important first step
towards establishing asymptotic regimes. Focusing on a ``critical''
fluid model (meaning, at least one link's capacity is saturated,
i.e., equal to the offered traffic load), the paper introduces the
notion of an {\it invariant state} --- if the fluid model starts in
such a state it will always stay in that state --- and shows that
under the law-of-large-numbers (fluid) scaling, the
bandwidth-sharing network operating under a general
proportional-fair allocation scheme converges to a fluid limit
(which, quite naturally, evolves over time to an invariant state).

In the broader framework of {processing networks} of Harrison
\cite{Harrison00,harrison},
 concurrent resource occupancy is one of
the highlighted features, and the main approaches are built around
Brownian network models.
Our study reported here is mostly motivated by the recent works of
Mandelbaum and Stolyar \cite{ms} and Stolyar
\cite{stolyar1,stolyar2}. In \cite{ms}, multi-class jobs are
processed by a set of ``flexible severs'', meaning each server can
process any job class albeit with different efficiency. It is shown
that the so-called generalized $c\mu$-rule
--- serving jobs according to a priority scheme that is determined
by the product of two factors: cost and service rate --- is
asymptotically optimal under a certain heavy-traffic regime. Similar
kinds of asymptotic optimality results are established in
\cite{stolyar1}, for a ``max-weight'' resource pooling scheme in a
general switch (e.g., including the crossbar switch as a special
case); and in \cite{stolyar2}, for a gradient scheduling algorithm
(which is more general than the proportional fair protocol), applied
also to a general switch.

\section{Utility Maximization and Cost Minimization}
\label{sec-model}

\subsection{The Network Model} \label{subsec-model}

Consider a network that consists of
a set of links, denoted $\call$.
There is also a set of routes, denoted $\calr$.
Each route $r\in\calr$ is a subset of links.
Denote $\ell\in r$ if link
$\ell$ is part of route $r$.

On each route, jobs arrive at the network
following a renewal process, independent of job arrivals on
all other routes.
Each job requires the simultaneous occupancy of all links
$\ell\in r$, for a period of service (or ``connection'')
time that is independent of all other jobs.
For each route $r$, denote the interarrival times
between consecutive jobs
as $u_r(i)$, and denote the amount of work (service requirement)
each job brings to the network
 as $v_r(i)$, $i=1,2,\dots$. Assume that
$u_r(i)$ ($i \ge 2$) are i.i.d.~random variables with mean
$1/\lambda_r$ and variance $a_r^2$, and that $v_r(i)$
are i.i.d.~random variables with mean $\nu_r$ and variance $b_r^2$.
 Denote the offered traffic load as $\rho := (\rho_r)_{r\in\calr}$, with
 $$\rho_r := \lambda_r \nu_r . $$
Assume $\lambda_r > 0$ and $\nu_r
> 0$ (hence, $\rho_r > 0$) for all $r\in\calr$.

\begin{figure}[htbp]
  \begin{center}
  \scalebox{0.3} {\includegraphics[angle=270]{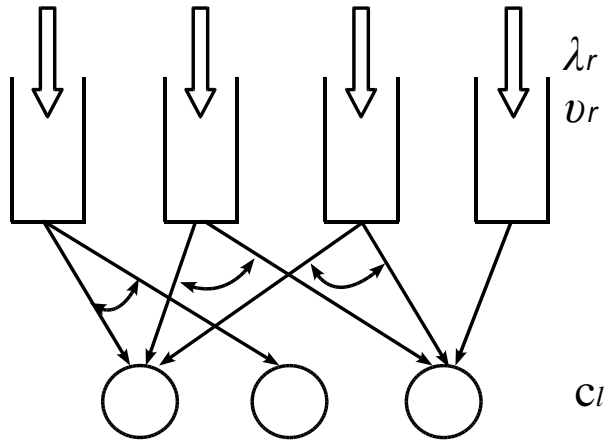}}
  \caption{A network example}
  \label{fig-net-example}
  \end{center}
  \end{figure}

An alternative view of the above network is depicted in Figure
\ref{fig-net-example}: there is
 a set of servers, corresponding to the links
$\call$; and a set of job classes, corresponding to the
routes $\calr$. To be processed, each
class-$r$ job requires the simultaneous occupancy of all the
servers $\ell\in r$.
Figure \ref{fig-net-example} shows
an example with four job classes and three
servers, with each of the first three classes
requiring two concurrent servers, while the fourth
class requires a single server.

A typical state of the network is denoted $\bn:=(n_r)_{r\in\calr}$,
where $n_r$ denotes the total number of class-$r$ jobs that are
present in the network. Each server $\ell\in\call$ has a given
capacity, $c_\ell$, which is shared among all job classes
$r\ni\ell$. More precisely, one job (if any) from each class
$r\ni\ell$ is processed at any time, while other jobs in the same
class waiting in a buffer.

The allocation of the service capacities takes place in each state,
denoted
$\G(\bn):=(\G_r(\bn))_{r\in\calr}$, where $\G_r(\bn)$ is the
capacity allocated to class $r$ when the network state is $\bn$. The
actual time needed to complete a job then depends on its service
requirement and the capacity allocated to it. Specifically, for the
$i$-th class-$r$ job mentioned above, provided it is being processed
in state $\bn$, then the amount of work $v_r(i)$ associated with it
is depleted at rate
$\G_r(\bn)$.

Let $M$ denote the set of all feasible allocations. Then, clearly
(omitting the argument $\bn$ from $\Lambda$),
\begin{eqnarray}
 \label{eq-lambda-fea1}
 M=\{\G:\; \sum_{r\ni\ell} \Lambda_r \le c_\ell,
 \,\forall\ell\in\call; \; \Lambda_r \ge 0 , \,\forall r\in\calr\} .
\end{eqnarray}

Throughout the paper, we shall assume that $\rho\in M$, i.e.,
there is enough capacity in the network to handle the offered load.
This is nothing more than the usual traffic condition,
often a necessary condition for stability.

\subsection{The Utility-Maximizing Allocation}
\label{sec-umax}

The utility-maximizing allocation refers to
the solution to the following optimization problem:
\begin{eqnarray}
\label{eq-u-util}
\max_{\G\in M}\sum_{r\in\calr} U_r(n_r,\Lambda_r) ,
\end{eqnarray}
where $M$, the set of all feasible allocations, follows the
specifications in (\ref{eq-lambda-fea1}); and $U_r(n_r,\Lambda_r)$,
$r\in\calr$, are utility functions defined on the two-dimension
nonnegative orthant  $\Re_+^2 = \{ x\in \Re^2: x_1, x_2\ge0 \}$. It
is standard to assume (e.g., \cite{ye}) that the utility functions
 are second-order differentiable
and satisfy the following conditions:
\begin{eqnarray}
 && U_r(0,\Lambda_r) \equiv 0 ; \label{eq-assume-0} \\
 && U_r (n_r,\Lambda_r)
  \mbox{ is strictly increasing and concave in } \Lambda_r ,
  \mbox{~~for } n_r > 0 ; \label{eq-assume-3} \\
 && \partial_2 U_r (n_r, \Lambda_r)
  \mbox{ in strictly increasing in } n_r \ge0 , \nonumber \\
  && \mbox{~~~~with } \partial_2 U_r (0, \Lambda_r)=0
    \mbox{ and } \lim_{n_r\rightarrow \infty}
       \partial_2 U_r (n_r, \Lambda_r)= +\infty,
  \mbox{~~for } \Lambda_r > 0 . \label{eq-assume-2}
\end{eqnarray}
 Here the operator $\partial_k f$ is a shorthand for the partial
derivative of function $f(\cdot)$ with respect to its $k$-th
variable.
Furthermore, we need another condition: An allocation is said to
satisfy the {\it partial radial homogeneity} (or {\it radial
homogeneity} for short) if for any scalar $a>0$, we have
\begin{eqnarray}
 && \Lambda_r (a\bn) = \Lambda_r (\bn) , \mbox{~~for any } r\in\calr
 \mbox{ with } n_r>0  . \label{eq-assume-4}
\end{eqnarray}
(Note, the qualifier, ``partial'', alludes to the fact that the
above is only required for each $r$ with $n_r>0$.)

Listed below are some examples of the utility function
widely used in modeling internet protocals:
\begin{eqnarray*}
 && U_r(n_r,\Lambda_r) = \beta_r n_r \log(\Lambda_r) , \nonumber\\
 && U_r(n_r,\Lambda_r) = - \beta_r \frac{n_r^2}{\Lambda_r} ,
  \mbox{ and } \nonumber\\
 && U_r(n_r,\Lambda_r) = \beta_r n_r^{\alpha}
 \frac{\Lambda_r^{1-\alpha}}{1-\alpha} ; \nonumber
\end{eqnarray*}
where $\alpha>0$
 and
$\beta_r>0$ are given parameters.
To motivate, consider the first utility function.
The optimal solution takes the following form
(from the first-order optimality equation):
$$\G_r=\frac{\beta_r n_r}{\eta_\ell}, \qquad r\in\ell,$$
where $\eta_\ell$ is the shadow price (Lagrangian multiplier)
associated with link $\ell$. That is, the optimal allocation among
job classes is proportional to the number of jobs present in the
network from each class. (Hence, this utility function is called
``proportionally fair''.) Note that the second utility function is a
special case of the third one when $\alpha=2$; and so is the first
one, in the sense that its maximizer coincides with the maximizer of
the third utility function when $\alpha\to 1$ (refer to \cite{bm}).
Also note that all three utility functions satisfy the conditions in
(\ref{eq-assume-0})-(\ref{eq-assume-2}); and the corresponding
optimal allocations also satisfy the condition in
(\ref{eq-assume-4}).

\medskip\noindent
{\bf Remarks.}
The familiar $c\mu$-rule can also be represented as a utility
maximizing allocation of a single server's (unit) capacity  among a
set of job classes, with a linear objective:
$$\sum_r U_r(n_r,\G_r) = \sum_r C_r \mu_r \ind(n_r>0)\G_r , $$
and a single constraint $\sum_r \G_r\le 1$. Hence, the solution is
to let $\G_r =1$ for the class $r$ (with at least one job present)
that corresponds to the largest $C_r\mu_r$ value.
 In the case of the generalized $c\mu$-rule where the linear cost
$C_r n_r$ is replaced by a general function $C_r(n_r)$ (see, e.g.,
\cite{ms,stolyar1}), the allocation
is the solution to the utility-maximizing problem $\max
\sum_{r} C_r'(n_r)\mu_r \Lambda_r$. The utility is then
$U_r(n_r, \Lambda_r) = C_r'(n_r)\mu_r  \Lambda_r$.
Note, however, in both cases the utility functions are linear in
$\G_r$; hence, they do not satisfy the {\it strict}
 concavity condition in (\ref{eq-assume-3}).
Furthermore, the strict increasingness (in $n_r$) of $\partial_2 U$
in (\ref{eq-assume-2}) implies strict convexity of
the cost function $C_r(n_r)$, and hence cannot be satisfied by
the $c\mu$-rule.
Consequently, it has been known (\cite{ye}) that the plain $c\mu$-rule
may not even guarantee stability, let alone asymptotic
optimality, in the kind of networks that we consider here.
On the other hand, for the generalized $c\mu$-rule, provided the
cost function $C_r(n_r)$ is strictly
convex and increasing, with $C_r(0)=0$, so that
the condition in (\ref{eq-assume-2}) is satisfied,
then the strict concavity (in $\G_r$) of the utility function,
the requirement in (\ref{eq-assume-3}), can be
 relaxed to non-strict concavity, and our main results below
will still hold.
 The proofs, however, will
involve more tedious, but non-essential, technical modifications.
For instance, the solution to the utility maximization problem in
(\ref{eq-u-util}) may not be unique; and hence, the allocation $\Lambda(n)$ can
be any one of the
 optimal solutions.
The radial homogeneity in (\ref{eq-assume-4}) needs to be modified
accordingly. The proofs of  the uniform attraction theorem and Lemma
\ref{lem-y-n} need to be modified as well, since Lemma \ref{lem-lambda-pty}
no longer holds.


\medskip

Returning to the utility maximization problem in (\ref{eq-u-util}),
we have the following optimality equations (which are both necessary
and sufficient, due to the concavity of the objective function
as assumed in (\ref{eq-assume-3})):
\begin{eqnarray}
 \label{uopteqn}
 \partial_2 U_r (n_r, \Lambda_r)=\sum_{\ell\in r}\eta_\ell,
 \quad \Lambda_r>0, r\in \calr; \qquad
 \eta_\ell\cdot(\sum_{r\ni\ell}\G_r -c_\ell)=0,
 \quad \ell\in \call ;
\end{eqnarray}
where $\eta_\ell \ge0$, $\ell\in \call$, are the Lagrangian
multipliers.

Furthermore, the following (partial) continuity property
of the utility-maximizing allocation $\Lambda(n)$, is known. It will
be used in the proofs of the uniform attraction theorem and Lemma
\ref{lem-y-n} below.
\begin{lem}
\label{lem-lambda-pty}
{\rm (\cite{ye}) Suppose a sequence of states $\{\bn^j, j =
1,2,\cdots \}$ converges to $\bn $ as $j \rightarrow \infty$. Then,
the solution to the utility maximization problem in
(\ref{eq-u-util}) converges accordingly:
\begin{eqnarray}
  &&\Lambda_r(n^j) \rightarrow \Lambda_r(n), \mbox{~~as } j
  \rightarrow \infty , \nonumber 
\end{eqnarray}
for each $r\in \calr$ such that $n_r>0$.
} 
\end{lem}

\subsection{A Cost Minimization Problem}
\label{sec-costmin}

For each $r\in\calr$, associated with the utility function $U_r$,
we define a cost function $C_r$
as follows:
\begin{eqnarray}
C_r(n_r,\G_r)= \nu_r\int_0^{n_r}  \partial_2 U_r(q, \G_r) d q ;
 \label{costutility} 
\end{eqnarray}
or, in differentiation form,
\begin{eqnarray}
\partial_1 C_r(n_r,\G_r)= \nu_r \partial_2 U_r(n_r, \G_r) .
\nonumber
\end{eqnarray}
(Recall, the utility functions are defined on $\Re_+^2$; hence, so
are the cost functions.)

For example, in the case of the third utility function above, we
have
 $$C_r(n_r,\G_r) = \frac{\beta_r \nu_r}{(1+\alpha) \G_r^{\alpha}}
 n_r^{1+\alpha}.$$
For the generalized $c\mu$-rule, where we have $U_r(n_r,
\G_r)=C_r'(n_r)\mu_r \G_r$, the relation in (\ref{costutility}) also
holds.

We say that the {\it heavy-traffic condition} holds, if the
offered-load vector $\rho$ is a maximal element of $M$, i.e., $\rho$
makes
a non-empty subset of constraints in $M$ binding:
\begin{eqnarray}
\label{htcond}
\sum_{r\ni\ell}  \rho_r = c_\ell, \quad \ell \in \call^* ;\qquad
 \sum_{r\ni\ell}  \rho_r < c_\ell, \quad \ell \not\in\call^* ;
\end{eqnarray}
for some $\call^* \subseteq \call$, and $\call^*\neq\emptyset$.
Intuitively, we shall refer to each link in $\call^*$ as a {\it
bottleneck} link.

Under the heavy-traffic condition, in particular, given the set of
bottleneck links $\call^*$, along with a given set of parameters
$w_\ell \ge 0$, $ \ell\in\call^*$,
 consider the following optimization problem:
\begin{eqnarray}
 \min_\bn && \sum_{r\in \calr} C_r(n_r,\G_r)
  \label{eq-workload-min} \\
 {\rm s.t.} && \sum_{r\ni\ell} \nu_r n_r \ge w_\ell, \qquad
 \ell \in \call^* , \nonumber\\
 && n_r \ge 0, \qquad
 r\in \calr .\nonumber
\end{eqnarray}
That is, for any given allocation $\Lambda$, we want to identify a
state $\bn$ under which the total cost over all routes is minimized
and the (average) workload at each bottleneck link $\ell\in\call^*$
is set to be greater than or equal to the required level $w_\ell$.

Let $\calr^*$ denote the set of routes $r$ such that $r\ni\ell$ for
some $\ell\in\call^*$, i.e., each route in $\calr^*$ involves at least one
bottleneck link. Then, clearly, for any route $r\not\in\calr^*$, we
must set $n_r=0$ in the optimal
solution to the minimization problem in (\ref{eq-workload-min}),
 since the cost $C_r$ is increasing in $n_r$,
following (\ref{eq-assume-0}) and (\ref{eq-assume-3}), along with
(\ref{costutility}). The remaining components of the optimal
solution, $(n_r,\,r\in\calr^*)$, can be obtained from the following
equations:
\begin{eqnarray}
\label{copteqn}
 \partial_1 C_r (n_r, \Lambda_r)=\nu_r\sum_{\ell\in r}\theta_\ell,
 \quad n_r>0, r\in \calr^*; \qquad
  \theta_\ell( \sum_{r\ni\ell} \nu_r n_r - w_\ell) = 0 ,
 \quad \ell\in \call^* ;
\end{eqnarray}
where $\theta_\ell$, $\ell\in \call^*$, are the Lagrangian
multipliers.


\begin{pro}
\label{pro:utilcost} {\rm Suppose the heavy-traffic condition in
(\ref{htcond}) holds; and suppose $\calr^*=\calr$, i.e., every route
involves at least one bottleneck link.
\begin{itemize}
\item[(i)] If $\bn^*$ is the optimal solution to
the cost minimization problem in (\ref{eq-workload-min}), with
$(\G_r=\rho_r)_{r\in\calr}$ in the cost function;
then, $(\G_r^*=\rho_r)_{r\in\calr}$  must be
the optimal solution to the utility maximization problem
in (\ref{eq-u-util}) with $\bn=\bn^*$ in the utility function.
\item[(ii)] Conversely, if $(\G^*_r=\rho_r)_{r\in\calr}$
is the optimal solution to the utility maximization problem;
then, the optimal solution to the cost minimization problem, with
$(\G_r=\G^*_r)_{r\in\calr}$ in the cost function,
must be $\bn^*=\bn$, with $\bn$
being the state vector in the
utility function, and with
$w_\ell=\sum_{r\ni\ell} \nu_r n_r$ for all $\ell\in\call^*$.
\end{itemize}
When $\calr^*\subset\calr$, both statements above still hold, with
the correspondence between the two optimal solutions,
$\G_r^*=\rho_r$ and $n_r^*=n_r$, restricted to $r\in\calr^*$.
}
\end{pro}

{\startb Proof.} First, consider $\calr^*=\calr$. Suppose the
condition in (i) is true. From the heavy-traffic condition in
(\ref{htcond}), we know $(\G_r=\rho_r)_{r\in\calr}$ is a feasible
solution to the utility maximization problem. To argue it is also
optimal, note that
 for each $r\in\calr$, we have,
$$\partial_2 U_r (n^*_r, \rho_r)=\nu^{-1}_r\partial_1 C_r (n^*_r,\rho_r)=
\sum_{\ell\in r}\theta_\ell,$$ where the first equation is due to
(\ref{costutility}), and the second one follows from (\ref{copteqn})
since $\bn^*$ is the optimal solution to the cost minimization
problem. Hence, letting
 $\eta_\ell=\theta_\ell$ for $\ell\in\call^*$,
and $\eta_\ell=0$ for $\ell\not\in\call^*$, will satisfy the
optimality equations in (\ref{uopteqn}).

To go the other way, follow the same argument. In particular, the
first set of equations in (\ref{uopteqn}) translates into the first
set of equations in (\ref{copteqn}), with $\theta_\ell=\eta_\ell$,
for $\ell\in\call^*$. The second set of equations in (\ref{copteqn})
holds automatically.

When $\calr^*\subset\calr$, the above arguments will only
apply to $r\in\calr^*$. Specifically, by setting
 $\eta_\ell=\theta_\ell$ for $\ell\in\call^*$,
we can conclude
that $(\G_r^*=\rho_r)_{r\in\calr^*}$ is (part of) the
optimal allocation to the utility maximization problem.
(For $r\not\in\calr^*$, $\G_r$ is not specified; and neither
is $\eta_\ell$ for $\ell\not\in\call^*$.)
Similarly, for the cost minimization problem, we have
$n^*_r=n_r$ for $r\in\calr^*$. For $r\not\in\calr^*$,
as appointed out earlier, it is optimal to have $n_r=0$.
\hfill $\qed$


\medskip

In the literature, the cost minimizer $n^*$ (of the cost
minimization problem with $\Lambda$ replaced by $\rho$) is often
referred to as a {\it fixed point} (\cite{ms,stolyar1}, or an {\it
invariant state} (\cite{kw}), provided the utility-maximizing
allocation in that state is $(\G^*_r=\rho_r)_{r\in\calr^*}$.
 We denote the fixed point corresponding to the
workload $w$ as $n^*(w)$.
 For convenience, we will suppress the second variable from the cost
function to write $C_r(n_r)$ for $C_r(n_r,\rho_r)$ and $C'_r(n_r)$
for $\partial_1 C_r(n_r,\rho_r)$, when $\Lambda_r = \rho_r$.

To close this subsection, we present the following lemma
on the continuity of the the fixed point $n^*(w)$, which will be
used in the proof of the main results (Theorems \ref{thm-fluid-ua}
and \ref{thm-opt-diff}) later.

\begin{lem}  \label{lem-fixPt-work-continue}
{\rm
The fixed point $\bn^*(w)$ as a function of the workload
$w=(w_\ell)_{\ell\in\call^*}$ is continuous in $w$.
}
\end{lem}

{\startb Proof.}
It suffices to show that for any sequence of workloads (vectors)
$\{w^{(j)}, j=1,2,...\}$ satisfying
\begin{eqnarray*}
 w^{(j)} \to w^{(0)} ,
 \qquad n^*(w^{(j)}) \to n' , \qquad
{\rm as}\quad j \to\infty , \nonumber
\end{eqnarray*}
we have $n' = n^*(w^{(0)})$.

First, the optimal value of the cost minimization problem in
(\ref{eq-workload-min}), $\sum_{r\in \calr} C_r(n_r^*(w)) $, as a
function of the workload $w$, is continuous in $w$. Therefore, we
have
\begin{eqnarray}
 && \sum_{r\in \calr} C_r(n_r^*(w^{(j)}))  \to
    \sum_{r\in \calr} C_r(n_r^*(w^{(0)}))
    \mbox{~~~~as } \; j \to  \infty .
    \nonumber
\end{eqnarray}
On the other hand, the cost function $C_r(n_r)$ is continuous in the
state $n_r$, $r\in {\cal R}$; and hence,
\begin{eqnarray}
  \sum_{r\in \calr} C_r(n_r^*(w^{(j)}))  \to
    \sum_{r\in \calr} C_r(n_r')
    \mbox{~~~~as } \; j \to  \infty .
    \nonumber
\end{eqnarray}
Therefore,
\begin{eqnarray}
 && \sum_{r\in \calr} C_r(n_r^*(w^{(0)}))
    = \sum_{r\in \calr} C_r(n_r') . \label{eq-2-3-30}
\end{eqnarray}
Next, for any real number $\delta >0$, let $S^\delta$ denote the set
of all states satisfying the following constraints,
\begin{eqnarray}
 n_r \ge 0, \quad r\in \calr;
 \qquad
 \sum_{r\ni\ell} \nu_r n_r  \ge  w_\ell^{(0)} - \delta,
    \quad \ell \in \call^* . \nonumber
\end{eqnarray}
 Then, it is directly verified that we have $n' \in S^\delta$.
Therefore, we
have
 $$n' \in \lim_{\delta \rightarrow 0} S^\delta = \cap_{\delta >0}
  S^\delta.$$
In other words, the state $n'$ satisfies the following constraints,
\begin{eqnarray}
 n_r \ge 0, \quad r\in \calr;
\qquad \sum_{r\ni\ell} \nu_r n_r \ge w_\ell^{(0)},
  \quad \ell \in \call^* . \nonumber
\end{eqnarray}
Finally, the above constraints, combined with the identity
(\ref{eq-2-3-30}), imply that the state $n'$ is also the optimal
solution to the cost minimization problem in (\ref{eq-workload-min})
with  $w$ replaced by $w^{(0)}$. Since the optimal solution is
unique, we must have $n'=n^*(w^{(0)})$.
\hfill $\qed$

\section{Fluid Limit and Uniform Attraction}
\label{sec-fluid}


Here we start with presenting the main processes associated with
the stochastic network introduced in the last section.
First, the two primitive processes are the renewal (counting)
processes associated with the job arrivals and the work
(service requirements) they bring into the network:
 $E(t)=(E_r(t))_{r\in\calr}$
and $S(y)=(S_r(y))_{r\in\calr}$, where
\begin{eqnarray}
E_r(t) &:=& \max\{i: u_r(1)+ \cdots +u_r(i) \le t\} ,
\label{eq-E-def}\\
S_r(y)&:=&\max\{i:v_r(1)+\cdots+v_r(i) \le y \} .
\label{eq-output-HOLPS}
\end{eqnarray}
The two derived processes that characterize, along with the
two primitive processes, the dynamics
of the stochastic network are the queue-length
process and the capacity allocation process:
$N(t)=(N_r(t))_{r\in\calr}$ and
$D(t)=(D_r(t))_{r\in\calr}$, where
\begin{eqnarray}
 N_r(t) &:=& N_r(0) + E_r(t) - S_r(D_r(t)) , \label{eq-N-E-S}\\
D_r(t) &:=& \int_0^t \Lambda_r(N(s)) 1_{ \{ N_r(s)>0 \} } ds.
\label{eq-D-Lambda}
\end{eqnarray}
Note that $N_r(0)$ is the initial
number of class-$r$ jobs in the system, $r\in\calr$.

In addition, we define two more derived processes:
$W(t)=(W_\ell(t))_{\ell\in\call^*}$
and $Y(t)=(Y_\ell(t))_{\ell\in\call^*}$, where
\begin{eqnarray}
 W_\ell(t) := \sum_{r\ni \ell} \nu_r N_r(t),
\nd
Y_\ell(t) := \sum_{r\ni \ell} (\rho_r t - D_r(t)) .
 \label{eq-2-dyn-y}
\end{eqnarray}
Clearly, $W(t)$ translates the queue-lengths into the (expected)
workloads, and $Y(t)$ keeps track of the unused capacities at the
bottleneck links. (Recall, as we noted in the last section, in a
fixed-point state, the utility-maximizing allocation at each
bottleneck link is $\G_r=\rho_r$.)

The following property of $Y(t)$ will be used later:
\begin{eqnarray}
 Y_\ell(t) \mbox{ is non-decreasing in }\, t\ge0,
 \mbox{~~~~for all } \ell\in L^* . \label{eq-2-700}
\end{eqnarray}
The above holds because from
(\ref{eq-D-Lambda}) and (\ref{eq-2-dyn-y}), we have
\begin{eqnarray}
 && Y_\ell(t)= \int_0^t
 \sum_{r\ni \ell} [ \rho_r- \Lambda_r(N(s)) 1_{\{N_r(s)>0\}} ] ds ;
 \nonumber
\end{eqnarray}
and
\begin{eqnarray}
 \sum_{r\ni \ell} \Lambda_r(N(s))\le c_\ell
 = \sum_{r\ni \ell} \rho_r ,
 \qquad \mbox{for } s\ge0, \ell\in\call^* .
 \nonumber
\end{eqnarray}

To describe the fluid limit and its uniform attraction property, we
introduce a sequence of networks,
 indexed by $k$.
Each of the networks is like the one introduced in the last
section,
but may differ in their arrival rates and mean service times
(which are also indexed by $k$).  We assume, for each $r\in\calr$,
\begin{eqnarray}
\label{flurates}
 \lambda_r^k \to \lambda_r  \nd
 \nu_r^k \to \nu_r,
\qquad{\rm as }  k\to\infty ;
\end{eqnarray}
and consequently,
$ \rho_r^k \to \rho_r$.
In addition, we assume
\begin{eqnarray}
 \frac{u_r^k(1)}{k} \rightarrow 0 \mbox{~~and~~}
 \frac{v_r^k(1)}{k} \rightarrow 0 ,
  \mbox{~~~~as } k \rightarrow \infty ;
 \label{eq-3-ini-cond}
\end{eqnarray}
i.e., under the fluid scaling the initial interarrival times and
initial residual work requirements are negligible.


We apply the standard fluid scaling to the primitive processes
associated with this sequence of networks:
$$\left( \bar E^k(t), \bar S^k(y)) \right):=
\left( \frac{1}{k} E^k( k t),   \frac{1}{k} S^k( k y)\right);$$
and similarly define the fluid-scaled version of the derived processes:
\begin{eqnarray*}
 \left(  \bar N^k(t), \bar D^k(t),
 \bar W^k(t), \bar Y^k(t) \right)
 =   \left(  \frac{1}{k} N^k( k t),\frac{1}{k} D^k( k t),
   \frac{1}{k} W^k( k t), \frac{1}{k} Y^k( k t) \right) .
\end{eqnarray*}

When $k\to\infty$ and under the assumption (\ref{eq-3-ini-cond}), we
know (e.g., \cite{cynet}, Chapter 5)
\begin{eqnarray} \label{primconv}
 \Big( \bar E^k(t), \bar S^k(y)\Big)
 \to (\g t, \mu y), \qquad {\rm u.o.c.},
\end{eqnarray}
where $\g:=(\g_r)_{r\in\calr}$ and $\mu:=(\nu_r^{-1})_{r\in\calr}$, and
the convergence is uniform on compact sets (u.o.c.) of $t\ge 0$.

The theorem below states that the sequence of derived processes also
approaches a limit, the fluid limit, which can be characterized as
follows.
Let
$\bar N (t):=(\bar N_r (t))_{r\in\calr}$ and
$\bar D (t):=(\bar D_r (t))_{r\in\calr}$, where
\begin{eqnarray}
 \bar N_r(t) &:=& \bar N_r(0) + \lambda_r t - \nu_r^{-1} \bar D_r(t) \ge 0 ,
 \label{eq-fluid-30} \\
  \bar D_r(t) &:=& \int_0^t \bar \Lambda_r(\bar N(s)) ds ;
 \label{eq-fluid-35}
\end{eqnarray}
and
\begin{eqnarray}
&&  \bar \Lambda_r(n) := \left\{
\begin{array}{ll}
   \Lambda_r(n)  & \mbox{if } n_r > 0,  \\
   \rho_r   & \mbox{if } n_r = 0 .
\end{array}
\right.\label{eq-fn-Lambda}
\end{eqnarray}
Furthermore, define
$\bar W(t):=(\bar W_\ell (t))_{\ell\in\call^*}$ and
$\bar Y(t):=(\bar Y_\ell (t))_{\ell\in\call^*}$,
where
\begin{eqnarray}
 \bar W_\ell(t) := \sum_{r \ni \ell} \nu_r \bar N_r(t)
\nd
\bar Y_\ell(t) := \sum_{r \ni \ell} ( \rho_r t - \bar D(t)) .
 \label{eq-fluid-50}
\end{eqnarray}
Then, clearly, we have (cf.\ (\ref{eq-2-700})):
\begin{eqnarray}
 && \bar Y_\ell(t) \mbox{ is non-decreasing in }\, t\ge0 . \label{eq-fluid-80}
\end{eqnarray}
Also note that the processes $\bar N(t), \bar D(t), \bar W(t)$ and $\bar Y(t))$ are
all Lipschitz continuous; and
\begin{eqnarray}
 && \bar W(t) = \bar W(0) + \bar Y(t)  \label{eq-fluid-70}
\end{eqnarray}
follows from
 (\ref{eq-fluid-30}) and
(\ref{eq-fluid-50}).


\begin{thm}
\label{thm-stoch2fluid}
{\rm {\bf (Fluid limit)}\, %
Given the utility-maximizing allocation (i.e., $\G (\cdot)$ is the
solution to (\ref{eq-u-util})), and suppose $\bar N^k(0)\to \bar
N(0)$.
Then, for any subsequence of these processes, there exists a further
subsequence, denoted ${\cal K}$, such that, along ${\cal K}$,
\begin{eqnarray}
 \left(  \bar N^k(t), \bar D^k(t),
 \bar W^k(t), \bar Y^k(t) \right) \to
 \left(  \bar N(t), \bar D(t),
 \bar W(t), \bar Y(t) \right)
\qquad {\rm u.o.c.}
\end{eqnarray}
with the fluid limit, $(\bar N(t),\bar D(t),\bar W(t),\bar Y(t))$,
following the specification in
(\ref{eq-fluid-30}), (\ref{eq-fluid-35}),
(\ref{eq-fn-Lambda}) and
(\ref{eq-fluid-50}).
} 
\end{thm}


\medskip

The proof of the above theorem is a slight variation of the proof of
Proposition 4.2 in Ye {\it et al.}~\cite{ye}, and thus is omitted.
We also remark that
the convergence of the two primitive processes in (\ref{primconv})
can be used as an alternative to the assumption
that these two are  renewal processes.


\medskip

The next theorem states that
 when $t\to\infty$, the fluid limit is
attracted to a fixed-point state, at which the cost
(the objective value in the cost minimization problem) is minimal
compared with any other fluid state with the same
workload.

We start with a definition and a lemma. Denote the cost objective
(with $\Lambda=\rho$) as
\begin{eqnarray}
 && \psi(n): = \sum_{r \in \calr} C_r(n_r). \label{psidef}
\end{eqnarray}
Below, we shall
focus on those differentiable points $t$ of the fluid limit, referred to
as ``regular'' times.
(Note that the fluid limit is Lipschitz continuous,
and hence differentiable almost everywhere.)

\begin{lem}\label{lem-fluid-ua}
{\rm Suppose $|\bar N(0)| \le B $ for some $B>0$. Then, $\bar N(t)$
is bounded by $\eta(B)$, where $\eta(\alpha)$ $(\alpha\ge0)$ is a
continuous increasing function  with $\eta(0)=0$;
and $ \bar W(t)$ is bounded and (component-wise) non-decreasing.
 Consequently, $\bar W(t)\to\bar W(\infty)$ as $t\to\infty$. }
\end{lem}

\noindent
{\bf Proof}.
The boundedness of $ \bar W(t)$ follows trivially from
the boundedness of $ \bar N(t)$: see (\ref{eq-fluid-50}).
In addition, from (\ref{eq-fluid-80})
and (\ref{eq-fluid-70}), we know $ \bar W(t)$ is non-decreasing in $t$.
Putting the two together, we have $\bar W(t)\to\bar W(\infty)$.
Hence, it suffices to the boundedness of $ \bar N(t)$.

We have, for any regular time $t\ge 0$,
\begin{eqnarray}
 \frac{d}{dt} \psi(\bar N(t)) & = & \sum_{r \in \calr} \dot {\bar N}_r(t)
  C'_r(\bar N_r(t)) \nonumber\\
 & = & \sum_{r \in \calr} \left( \rho_r -\bar \Lambda_r(\bar N(t)) \right) \cdot
 \partial_2 U_r(\bar N_r(t), \rho_r) \nonumber\\
 & = & \sum_{r \in \calr_t^+} \left( \rho_r - \Lambda_r(\bar N(t)) \right) \cdot
 \partial_2 U_r(\bar N_r(t), \rho_r)  \nonumber\\
 & \le & \sum_{r \in \calr_t^+} \left( U_r(\bar N_r(t), \rho_r)
   - U_r(\bar N_r(t), \Lambda_r(\bar N(t))) \right) \nonumber\\
 & = & \sum_{r \in \calr} \left( U_r(\bar N_r(t), \rho_r)
   - U_r(\bar N_r(t), \Lambda_r(\bar N(t))) \right) \le 0 ,
\label{eq-80}
\end{eqnarray}
where, the second equality follows from (\ref{eq-fluid-30}) and
(\ref{eq-fluid-35}); the third equality from (\ref{eq-fn-Lambda}),
with $\calr_t^+$ denoting the set of routes with a positive fluid
level at time $t$, i.e., $ \calr_t^+ = \{ r\in \calr: \bar N_r(t) >
0 \}$; the first inequality is due to concavity of the utility
function; and the second inequality is
due to the fact that $\Lambda(\bar N_r(t))$ is the optimal solution to the
utility maximization problem in  (\ref{eq-u-util}).
(Specifically, $( \Lambda_r(\bar N(t)) - \rho_r)_{r\in\calr}$
is the ascent direction of the objective function.)
Hence, $\psi(\bar N(t))\le \psi(\bar N(0))$ for all $t$. From
(\ref{psidef}), we can then conclude that $\bar N(t)$ is also
bounded, since $C_r (n_r)$, $r\in\calr$, are all strictly increasing
and unbounded functions.

Furthermore, define for any $\alpha\ge0$
\begin{eqnarray}
 && \tilde \eta(\alpha) = \max_{|n| \le \alpha} \psi(n)
  \mbox{~~and} \nonumber\\
 && \eta(\alpha) =
 \sum_{r\in {\cal R}} C_r^{-1}(\tilde \eta(\alpha)) . \nonumber
\end{eqnarray}
As the cost $C_r$ is strictly increasing, its inverse used above is
well defined and is strictly increasing. It is direct to check that
the function $\eta$ defined above is increasing with $\eta(0) = 0$.
 Now, we have for all $t\ge0$
\begin{eqnarray}
 && \psi(\bar N(t)) \le \psi(\bar N(0)) \le \tilde \eta(B), \nonumber
\end{eqnarray}
which implies (cf. (\ref{psidef})), for all $r\in {\cal R}$,
\begin{eqnarray}
 && C_r(\bar N_r(t)) \le \tilde \eta(B), \nonumber \\
 && \bar N_r(t) \le C_r^{-1}(\tilde \eta(B)), \nonumber
\end{eqnarray}
and finally
\begin{eqnarray}
 && |\bar N(t)| \le \sum_{r\in {\cal R}} C_r^{-1}(\tilde \eta(B))
    = \eta(B). \nonumber
\end{eqnarray}
 \hfill $\qed$

\medskip

Recall, $n^*(w)$ denotes the optimal solution to the cost
minimization problem in (\ref{eq-workload-min}). Given a state $\bar
N(t)$, the corresponding workload is $\bar W(t)$. If we set $w=\bar
W(t)$ in the cost minimization problem, we can write it more
explicitly as $w(\bar N(t))$. Then, the minimizer of the cost
objective can be expressed as $n^*(w(\bar N(t))$ and the
corresponding objective value as $\psi(n^*(w(\bar N(t)))$. To
lighten up notation, we shall write this simply as $\psi^*(\bar
N(t))$.
Note the difference between
$\psi^*(\bar N(t))$ and
 $\psi (\bar N(t))$: the latter is simply the cost objective
evaluated at a give $n=(\bar N_r(t))_{r\in\calr}$, i.e.,
the objective value of a mere {\it feasible} solution to
(\ref{eq-workload-min}); whereas $\psi^*(\bar N(t))$ is the
objective value of the optimal solution. Both, however, correspond
to the same required workload $w=w(\bar N(t)$.

Note that since $\bar W(t)$ is non-decreasing in $t$,
as we have shown in the last lemma, and
$C_r$'s are increasing functions, clearly
$\psi^*(\bar N(t))$ is also non-decreasing in $t$.
This is in contrast to the non-increasing property of
$\psi(\bar N(t))$ as established in
(\ref{eq-80}).
Indeed, the key to the proof of the theorem below is the
following Lyapunov function:
\begin{eqnarray}
\label{lyap}
L(n)=\psi (n)-\psi^*(n),
\end{eqnarray}
where, following the above notation $\psi^*(n)=\psi (n^*(w(n))$.

\begin{thm}
\label{thm-fluid-ua}
{\rm
{\bf (Uniform attraction)}\,
Suppose
$|\bar N(0)| \le B $ for some constant $B>0$.
Then, there exist a time,
$T_{B,\epsilon}$, and an
attraction state, $\bar N(\infty)$,
such that
\begin{eqnarray}
 && |\bar N(t) -\bar N(\infty)| < \epsilon ~~~~\mbox{for all }
 t > T_{B,\epsilon} \label{eq-fluid-attaction} .
\end{eqnarray}
Furthermore, the attraction state is a fixed point, i.e.,  $\bar
N(\infty)=n^*(w)$, with $w=\bar W(\infty)$ being the limit specified
in Lemma \ref{lem-fluid-ua}. }
\end{thm}

\noindent
{\bf Proof}.
Consider any given regular time $t\ge0$. We first show that
 for any $\delta>0$, there exists a $\sigma>0$ such that
\begin{eqnarray}
 |\bar N(t) - n^*(\bar W(t))| \ge \delta  \qquad\imp\qquad
 \frac{d}{dt}L(\bar N(t)) \le -\sigma, \label{eq-120}
\end{eqnarray}
where $L(n)$ is the Lyapunov function defined in (\ref{lyap}).
Since as argued above, $\psi(\bar N(t))$ is non-increasing
and $\psi^* (\bar N(t))$ non-decreasing,
it suffices to show the statement in (\ref{eq-120}) with
$\frac{d}{dt}L(\bar N(t))$ replaced by $\frac{d}{dt}\psi(\bar N(t))$.

From (\ref{eq-80}), we know $\frac{d}{dt}\psi(\bar N(t))\le 0$,
with equality holding if and only if $\bar N(t)$ is a fixed point,
i.e., $\bar N(t) = n^*(\bar W(t))$. This, along with the observation
that the function
 $$ f(n) : =
 \sum_{r \in {\cal R}_n} \left( \rho_r -\bar \Lambda_r(n) \right)
 \cdot \partial_2 U_r(n_r, \rho_r),
 \quad{\rm where}\quad {\cal R}_n := \{ r\in\calr: n_r > 0 \},$$
is continuous in $n$ (cf. Lemma \ref{lem-lambda-pty}), yields the
desired conclusion.

The statement in (\ref{eq-120}) implies that
for any $\epsilon > 0$, there
exists a time $T_{B,\epsilon}^1$ such that
\begin{eqnarray}
 |\bar N(t) - n^*(\bar W(t)) | < \frac{\epsilon}{2}
 \mbox{~~for all }\;  t> T_{B,\epsilon}^1. \label{eq-150a}
\end{eqnarray}
On the other hand, the convergence of $\bar W(t)\to\bar W(\infty)$
in Lemma \ref{lem-fluid-ua}, along with the
continuity of the fixed point in Lemma
\ref{lem-fixPt-work-continue}, guarantees that
\begin{eqnarray}
 | n^*(\bar W(t)) - n^*(\bar W(\infty)) | < \frac{\epsilon}{2}
  \mbox{~~for all } \; t> T_{B,\epsilon}^2. \label{eq-150}
\end{eqnarray}
Hence, the desired result in (\ref{eq-fluid-attaction})
 follows from
letting $T_{B,\epsilon} = \max \{ T_{B,\epsilon}^1,
T_{B,\epsilon}^2\} $, and
combining (\ref{eq-150a}) and (\ref{eq-150}).
 \hfill $\qed$

\begin{cor}
\label{cor-fluid-ua} {\rm (a) If $\bar N(0)$ is bounded ($\bar N(0)
\le B$), then there exists some $T(=T_{B,\epsilon})>0$, such that,
for any $\epsilon
>0$, the following holds for all $t\ge T$:
\begin{eqnarray}
 | \bar N(t) - n^* (\bar W(t))| < \epsilon . \nonumber
\end{eqnarray}
(b) If $\bar N(0)$ is a fixed point state,
then $\bar N(t)=\bar N(0)$ for
all $t\ge0$.
In particular, if $\bar N(0) = 0$,
then $\bar N(t)=0$ for all $t\ge0$.
}
\end{cor}

\noindent
{\bf Proof}.
Part (a) simply restates what's already proved in
(\ref{eq-150a}). Part (b) follows directly from the
property of the Lyapunov function proved above.
 \hfill $\qed$


\medskip\noindent
{\bf Remark.} Following Theorem \ref{thm-fluid-ua}, when all links
are non-bottlenecks, i.e., $\sum_{r\ni \ell} \rho_r < c_\ell$ for
all $\ell \in {\cal L}$, we have $\bar N(\infty) =0$. (In fact, we
have  $\bar N(t) =0$, $t\ge T$, for some $T>0$.) This recovers the
stability result of Ye {\it et al.}~\cite{ye}. Although the network
model in \cite{ye} allows complete processor sharing -- the link
capacities are shared among {\it all} jobs present in the network,
whereas ours is essentially a head-of-the-line processor-sharing
mechanism, the two models are equivalent when the service
requirements follow exponential distributions, which is indeed
assumed in \cite{ye}. Therefore, our results here are more general,
in allowing general distributions of the processing requirements

\section{Diffusion Limit and Heavy-Traffic Optimality}
\label{sec-diffusion}

As in Section \ref{sec-fluid}, we introduce a sequence of networks
indexed by $k$. In addition to the limits of the arrival rates and
mean processing times in (\ref{flurates}), we assume the existence
of the following limits as $k \to \infty$:
\begin{eqnarray}
 k(\lambda_r^k - \lambda_r) \to \theta_{r,\lambda}
 \nd
 k(\nu_r^k - \nu_r) \to \theta_{r,\nu} ,\qquad r\in\calr .
 \label{eq-4-theta}
\end{eqnarray}
Consequently, we have,
\begin{eqnarray}
 k (\rho_r^k - \rho_r ) =  k (\lambda_r^k \nu_r^k - \lambda_r \nu_r )
 \to \lambda_r \theta_{r,\nu} + \nu_r \theta_{r,\lambda} :=
 \theta_{r,\rho} . \label{eq-4-def-theta}
\end{eqnarray}
Moreover, we also need to
assume the existence of the limits of the standard deviations
 of the inter-arrival times and service requirements:
\begin{eqnarray}
 a_r^k \to a_r  \nd
b_r^k \to b_r , \qquad r\in\calr . \label{eq-4-40}
\end{eqnarray}

We apply the standard diffusion scaling (along with centering)
to the primitive processes, which results in the following:
\begin{eqnarray}
 \hat E_r^k(t) :=  \frac{1}{k} \left( E_r^k( k^2 t) - \lambda_r^k k^2 t \right),
 \qquad \hat S_r^k(t) :=  \frac{1}{k} \left( S_r^k( k^2 y) - (\nu_r^k)^{-1} k^2 y \right).
 \label{eq-diff-102}
\end{eqnarray}
Apply the same diffusion scaling to the derived processes:
\begin{eqnarray}
\hat N_r^k(t) :=  \frac{1}{k}  N_r^k( k^2 t) ,
\quad \hat W_\ell^k(t) :=  \frac{1}{k}  W_\ell^k( k^2 t) ,
\quad \hat Y_\ell^k(t) :=  \frac{1}{k}  Y_\ell^k( k^2 t) .
\label{eq-diff-108}
\end{eqnarray}

Next, we re-express
the unscaled workload process for the $k$-th network
as follows:
\begin{eqnarray*}
 W_\ell^k(t) &=& \sum_{r\ni\ell} \nu_r^k N_r^k(t) \\
 & = & \sum_{r\ni\ell} \nu_r^k N_r^k(0)
 + \sum_{r\ni\ell} \nu_r^k \left( [E_r^k(t)-\lambda_r^k t]
   - [ S_r^k(D_r^k(t)) - (\nu_r^k)^{-1} D_r^k(t)] \right) \\
 &&  + \sum_{r\ni\ell} (\rho_r^k -\rho_r) t
  +  \sum_{r\ni\ell} (\rho_r t - D_r^k(t)) .
\end{eqnarray*}
Applying the diffusion scaling to both sides of the
above equation, we have
\begin{eqnarray}
 \hat W_\ell^k(t) := \frac{1}{k}  W_\ell^k( k^2 t)
 = \hat X_\ell^k(t) + \hat Y_\ell^k(t) , \label{eq-diff-110}
\end{eqnarray}
where
\begin{eqnarray}
 \hat X_\ell^k(t) &:=& \frac{1}{k} \sum_{r\ni\ell} \nu_r^k N_r^k(0)
   + \sum_{r\ni\ell} \nu_r^k
   \left( \hat E_r^k(t) - \hat S_r^k({\tilde D}_r^k(t)) \right)
    + \sum_{r\ni\ell} k (\rho_r^k -\rho_r) t  , \label{eq-diff-130} \\
 \hat Y_\ell^k(t) &=&
\frac{1}{k} \sum_{r\ni\ell} [\rho_r k^2 t - k^2 {\tilde D}_r^k(t)]; \label{eq-diff-140}
\end{eqnarray}
and ${\tilde D}_r^k(t) := \frac{1}{k^2} D_r (k^2 t)$ is a variation
of the fluid-scaled process $\bar D_r(t)$.
  Note that, $\hat Y_\ell^k(t)$ in (\ref{eq-diff-140}) is consistent with
what it was originally defined in (\ref{eq-diff-108}). It is also
completely analogous to its fluid-scaled counterpart in
(\ref{eq-fluid-50}). 
 Indeed, we also have (cf.\ (\ref{eq-2-700})),
for each $k$,
\begin{eqnarray}
 && Y_\ell^k(t) \mbox{ is non-decreasing in}\, t\ge0,
 \mbox{~~~~for all } \ell\in L^*  . \label{eq-diff-146}
\end{eqnarray}

Following conventional arguments
(the functional central limit theorem, the random time-change
theorem, etc.; refer to \cite{cynet}, Chapter 5),
the processes $\hat E_r^k(t)$,
$\hat S_r^k({\tilde D}_r^k(t))$ and $\hat X_\ell^k(t)$, for
$r\in\calr$, converge
weakly to some Brownian motion.
Specifically, the limit of $\hat E^k_r(t)$, denoted
$\hat E_r(t)$, is a
Brownian motion with zero mean and variance $\lambda_r^3 a_r^2 $.
Assume, for ease of exposition, $N^k_r(0)=0$ for all $r\in\calr$ and
all $k$.
(More general initial conditions can be handled by addressing the
convergence of the sequence of initial states,
and the associated inter-arrival times  and  service
requirements of the initial jobs. Details are similar to those in \cite{ms}.)
Then, ${\tilde D}_r^k(t)\to\rho_r t$.
Hence,
$\hat S_r^k({\tilde D}_r^k(t))$ converges to
$\hat S^k(\rho_r t)$, which is a
Brownian motion with zero mean and variance $\nu_r^{-3} b_r^2 $; and
  $\hat X_\ell^k(t)$ converges to the following limit:
\begin{eqnarray}
\hat X_\ell(t) :=
    \sum_{r\ni \ell} \theta_{r,\rho} t  +
    \sum_{r\ni \ell} \nu_r \left( \hat E_r(t) - \hat S_r(\rho_r t)
    \right),
\label{eq-diff-uocW}
\end{eqnarray}
which is a Brownian motion with
drift $\sum_{r\ni
\ell} \theta_{r,\rho}$  and variance $ \sum_{r\ni \ell} (\nu_r^2
\lambda_r^3 a_r^2 + \nu_r^{-1} b_r^2 )$.


For the derived processes, denote their limits (to be proved in the
next theorem under a single-bottleneck condition) as follows:
\begin{eqnarray}
\hat N(t) = (\hat N_r(t))_{r\in \calr},\quad
  \hat W(t) = (\hat W_\ell(t))_{\ell\in \call^*}, \quad
  \hat Y(t) = (\hat Y_\ell(t))_{\ell\in \call^*}.  \nonumber
\end{eqnarray}
These processes are characterized by the following
relations, for all $t \ge 0$, $\ell\in\call^*$ and $r\in\calr$:
\begin{eqnarray}
 && \hat W_\ell(t) = \hat X_\ell(t) + \hat Y_\ell(t) \ge 0 ;
 \label{eq-RBM-10} \\
 && \hat Y_\ell(t) \; \mbox{is non-decreasing in}\, t;
\quad{\rm with}\; \hat Y_\ell(0)=0;
\label{eq-RBM-50} \\
&& \int_0^\infty  \hat W_\ell(t) d \hat Y_\ell(t) = 0 ; \label{eq-RBM-50a}\\
  && \hat N_r(t) = n_r^*(\hat W(t)) .
 ~~~~\mbox{ for } t\ge0 , r\in {\cal R} . \label{eq-RBM-fixed}
\end{eqnarray}
It is known (e.g., \cite{cynet}, Chapter 6) that given $\hat X (t)$
the relations in (\ref{eq-RBM-10}), (\ref{eq-RBM-50}) and
(\ref{eq-RBM-50a}), which constitute the so-called Skorohod problem,
uniquely define the processes $\hat Y (t)$ and $\hat W (t)$: $\hat Y
=\Phi(X)$ and $\hat W =\Psi(\hat X)$, with $\Phi (\cdot)$ and $\Psi
(\cdot)$ being Lipschitz continuous mappings. In particular, when
$\hat X (t)$ is a Brownian motion, $\hat W (t)$ is a {\it reflected}
Brownian motion (RBM), and $\hat Y (t)$ is the associated {\it
regulator}.

\begin{thm} \label{thm-opt-diff}
{\rm
Suppose the heavy-traffic condition in (\ref{htcond}) holds, with a
{\it single} bottleneck link, i.e., $L^* = \{ \ell^* \}$ is a
singleton set. Under the utility-maximizing allocation, we have the
following results.
\begin{itemize}
\item[(a)]
{\bf (Diffusion Limit)}\,
The following weak convergence holds when $k\to\infty$:
\begin{eqnarray}
 \left( \hat W^k(t) , \hat Y^k(t) , \hat N^k(t) ) \right)
 \Rightarrow \left( \hat W(t) , \hat Y(t) , \hat N(t) \right),
\nonumber
\end{eqnarray}
with the limit following the specifications in (\ref{eq-RBM-10})
through (\ref{eq-RBM-fixed}); in particular, $\hat W$ is a
single-dimensional RBM, $\hat Y$ is the associated regulator, and
$\hat N$ is the fixed-point to the cost minimization problem with
$w=\hat W$.

\item[(b)]
{\bf (Asymptotic Optimality)}\, $\hat W$ and $\hat N$ are
minimal in the following sense:
Let $\hat W^{k,G}$ and $\hat N^{k,G}$ denote the processes associated
with any feasible allocation scheme $G$. Then, for all $t\ge 0$, we have
\begin{eqnarray}
 \liminf_{k\to\infty} \ex[ \hat W^{k,G}(t)] \ge \ex [\hat W (t)]
\nd  \liminf_{k\to\infty} \sum_{r\in\calr}\ex[ C_r(\hat N_r^{k,G}(t))]
 \ge \sum_{r\in\calr}\ex[ C_r(\hat N_r(t))] .
\label{eq-4-55}
\end{eqnarray}
\end{itemize}
}
\end{thm}

\medskip

The proof of the above theorem is deferred to Section
\ref{sec-diff-proof1}. Below, we first motivate the
single-bottleneck condition, a key requirement in the theorem.

\subsection{The Single-Bottleneck Condition}
\label{subsec-pooling}

Here we show that the single-bottleneck condition
in Theorem \ref{thm-opt-diff}
 is equivalent to the
so-called resource pooling condition, which has been widely used in
related studies on the heavy-traffic limits in various stochastic
processing network; refer to, e.g.,
\cite{Harrison00,harrison,HarrisonLopez99,HarrisonVanMieghem97,KellyLaws93,Laws92,ms,stolyar1,Williams00},
among others.

Following Harrison \cite{Harrison00}, the resource pooling condition can
be stated as the uniqueness of the dual optimal solution to
what's known as a ``static planning LP'' (linear program), which
in our context takes the following form. First, the primal LP,
with $\xi$ and $(\G_r)_{r\in\calr}$ as decision
variables:
\begin{eqnarray*}
{\rm (P)}~~~~ \max && \xi \\
{\rm s.t.} && -\G_r+\rho_r\xi \le 0 , \quad  r\in\calr ; \\
 && \sum_{r\ni\ell}\G_r\le c_\ell,  \quad  \ell\in\call ; \\
&& \xi \ge 0, \quad \G_r\ge 0, \quad  r\in\calr .
\end{eqnarray*}
Next, the dual, with $p_r$, $r\in\calr$, and $\pi_\ell$,
$\ell\in\call$, as variables:
\begin{eqnarray*}
{\rm (D)}~~~~ \min && \sum_{\ell\in\call} c_\ell\pi_\ell \\
{\rm s.t.}
 && \sum_{r\in\calr} \rho_r p_r \ge 1  ; \\
&& -p_r+\sum_{\ell\in r}\pi_\ell \ge 0 , \quad  r\in\calr ; \\
&& p_r\ge 0, \quad  r\in\calr; \quad  \pi_\ell \ge 0, \quad
\ell\in\call.
\end{eqnarray*}
Here, the dual variables $\pi_\ell$, $\ell \in {\cal L}$, are
 the shadow prices of the link capacities; and $p_r$
is the marginal cost for processing a class-$r$ job, $r\in\calr$.
The resource pooling condition can then be stated as the uniqueness
of the $(p_r)_{r\in\calr}$ part of the dual optimal solution. (No
requirement on the $\pi$ part.)

For our LP presented above, we claim the primal optimal solution is:
\begin{eqnarray}
\label{psol} \xi=1 ; \quad \G_r =\rho_r, \quad  r\in\calr ;
\end{eqnarray}
whereas the dual optimal solution is characterized as follows:
\begin{eqnarray}
\label{dsol2} p_r=\sum_{\ell\in r}\pi_\ell= \sum_{\ell\in
r,\,\ell\in\call^*} \pi_\ell ,  \qquad r\in\calr ;
\end{eqnarray}
and
\begin{eqnarray}
\label{dsol1} \sum_{\ell\in\call^*} c_\ell\pi_\ell =1; \qquad
\pi_\ell=0, \quad \ell\not\in\call^* .
\end{eqnarray}
To justify, it is straightforward to
check that the above satisfies:
(a) primal feasibility, (b) dual feasibility, and (c)
complimentary slackness.
Indeed, the only less-than trivial part is to check
the dual constraint, $\sum_{r\in\calr} \rho_r p_r
\ge 1$. It follows from
\begin{eqnarray*}
\sum_{r\in\calr} \rho_r p_r= \sum_{r\in\calr} \rho_r \sum_{\ell\in
r}\pi_\ell =\sum_{\ell\in \call}\pi_\ell \sum_{r\ni\ell} \rho_r
=\sum_{\ell\in \call^*}\pi_\ell \sum_{r\ni\ell} \rho_r
=\sum_{\ell\in \call^*}\pi_\ell c_\ell = 1,
\end{eqnarray*}
where the second last equality follows from $\sum_{r\ni\ell}\rho_r = c_\ell$
for $\ell\in\call^*$.

From (\ref{dsol2}) and (\ref{dsol1}),
it is clear that when there is a single bottleneck
link, i.e., $\call^*=\{\ell^*\}$, then the following dual optimal
solution is unique:
\begin{eqnarray*}
\pi_{\ell^*}=c^{-1}_{\ell^*};
\qquad \pi_\ell=0, \quad \ell\neq\ell^* ;\\
p_r=\pi_{\ell^*}=c^{-1}_{\ell^*}, \quad r\ni\ell^* ; \qquad p_r=0,
\quad r\not\ni\ell^* .
\end{eqnarray*}
Otherwise, it is not, as we can choose different
$(\pi_\ell)_{\ell\in\call^*}$ values
to satisfy the
first equation in (\ref{dsol1}), and thereby lead to
different $p_r$ values following (\ref{dsol2}).

As pointed out by Stolyar (Theorem 3 of \cite{stolyar1}), the
uniqueness of the optimal dual solution is equivalent to the
uniqueness of a vector $(p_r)_{r\in\calr}$ that satisfies:
\begin{eqnarray}
  && \sum_{r\in R} p_r \rho_r
  = \max_{\Lambda \in M} \sum_{r\in R} p_r \Lambda_r .
  \nonumber
\end{eqnarray}
Clearly, ${p}:=(p_r)_{r\in\calr}$
is the normal vector to the outer boundary of $M$
at $\rho$.
Since the outer boundary is formed
by a set of hyperplanes, each of which corresponds to a link,
${p}$ is unique if and only if it lies in the {\it interior}
of a face of $M$. This provides a geometric interpretation
of the single-bottleneck condition.


\medskip

We conclude this part by the following lemma, which says that
  when the network state $n$ is
close to a fixed point, the capacity of the single
bottleneck link will be fully utilized.
 The proof of the lemma is postponed to the Appendix;
the lemma itself will be used in the next subsection.


\begin{lem}
\label{lem-y-n}
{\rm
Suppose $\ell^*$ is the only bottleneck link in the network.
Then, for any given constant $\epsilon
>0$, there exists another constant
$\sigma=\sigma(\epsilon) >0$, such that, for any state $n$ with
workload
$w_{\ell^*}=\sum_{r\ni\ell^*}
\nu_r n_r$ satisfying
\begin{eqnarray}
|n - n^*(w)| \le \sigma ,\nd
 w_{\ell^*} \ge \epsilon  ;
 \nonumber
\end{eqnarray}
we will have $\sum_{r\ni\ell^*} \Lambda_r(n) = c_{\ell^*}$.
}
\end{lem}

\subsection{A Key Lemma}
\label{sec-lem-key-diff}

To simplify notation, since there is a single bottleneck link
$\ell^*$, we will drop the subscript $\ell^*$ from $\hat
X_{\ell^*}^k(t)$, $\hat Y_{\ell^*}^k(t)$, and $\hat W_{\ell^*}^k
(t)$, and from their (proposed) limits, $\hat X_{\ell^*} (t)$, $\hat
Y_{\ell^*}(t)$, and $\hat W_{\ell^*} (t)$. For the rest of this
subsection and the next one, these are all scalar processes. In
addition, we find it more convenient to adopt a sample-path approach
based on the Skorohod representation theorem, i.e., to turn the weak
convergence into a probability one convergence (u.o.c.) of suitable
copies of the processes on a common probability space. Indeed, in
the rest of this section and the next section, we shall focus on a
given sample path.


  \begin{figure}[htbp]
  \begin{center}
  \scalebox{0.3} {\includegraphics[angle=270]{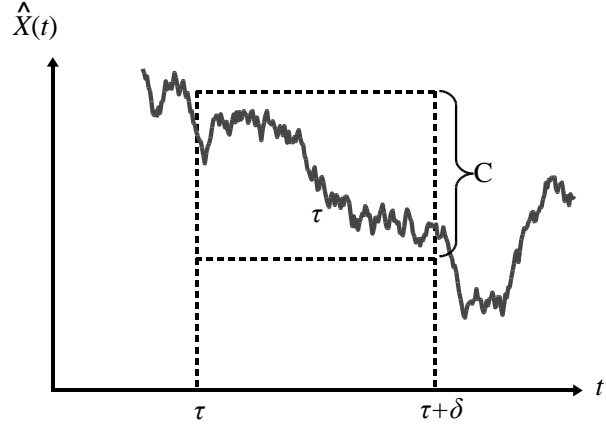}}
  \caption{Parameters $\tau$, $\delta$ and $C$}
  \label{fig-key-param}
  \end{center}
  \end{figure}


We shall focus on the time interval $[\tau , \tau+\delta]$, where
$\tau\ge0$ and $\delta >0 $. Let $T >0$ be a fixed time of a certain
magnitude to be specified later. Let the index $k$ be a large
integer. Divide the time interval $[\tau, \tau+\delta]$ into equal
segments of length $T/k$, a total of $k\delta / T$ such segments.
(Without loss of generality, assume $\delta / T$ is an integer, and
hence, so is $k\delta / T$.) Then, for any $t\in [\tau,
\tau+\delta]$, we can write it as $t=\tau+(jT+u)/k$ for some $j=0,
\cdots , k\delta / T$ and $u \in [0, T]$. Therefore, we write
\begin{eqnarray}
\hat W^k (t)= \hat W^k (\tau+\frac{jT+u}{k})
 = \frac{1}{k} W^k(k (k\tau + jT + u) )
 :=\bar W^{k,j}(u) , \qquad u\in[0,T], \; j\le \frac{k\delta}{T}.
  \label{eq-4-120}
\end{eqnarray}
That is, for each time point $t$, we will study the behavior of
$\hat W^k (t)$  through the fluid process,  $\bar W^{k,j}(u)$, over
the time interval $u\in [0,T]$. Similarly define $\bar N^{k,j}(u)$
and $\bar Y^{k,j}(u)$ as the fluid ``magnifiers'' of $\hat N^k (t)$
and $\hat Y^k (t)$.

\begin{figure}[htbp]
  \begin{center}
  \scalebox{0.3} {\includegraphics[angle=270]{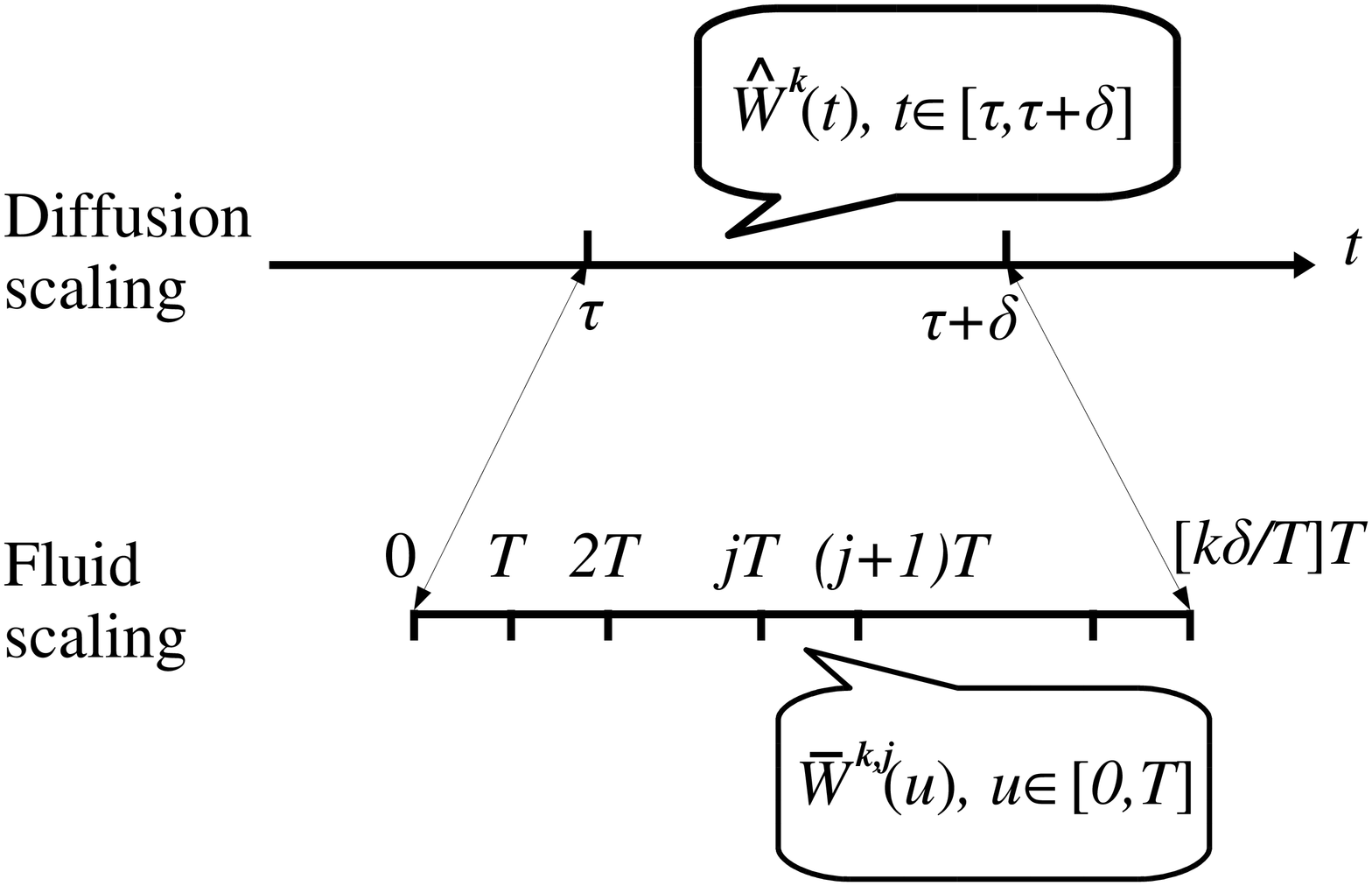}}
  \caption{Rescaling of processes}
  \label{fig-rescale}
  \end{center}
  \end{figure}

The above rescaling of $\hat W^k (t)$ is illustrated as an example
in Figure \ref{fig-rescale}. This rescaling technique enables one to
investigate the structure of diffusion scaled processes (e.g., $\hat
W^k(t)$) using the available results concerning the fluid scaled
processes (e.g., $\bar W^k(t)$). Such a technique appeared in
various forms in
\cite{Bramson98,ChenYe01,ChenZhang98b,ChenZhang00b,ms,stolyar1}.

\begin{lem} \label{lem-diff-key}
{\rm
Consider the time interval $[\tau , \tau+\delta]$, with $\tau\ge 0$
and $\delta >0$; pick a constant $C>0$ such that
\begin{eqnarray}
 \sup_{t', t'' \in [\tau , \tau+\delta]}
 | \hat X(t') - \hat X(t'') | \le C ; \label{eq-4-C-bdd}
\end{eqnarray}
and suppose
\begin{eqnarray}
\lim_{k\rightarrow\infty} \hat W^k(\tau) = \chi ,
\nd
\lim_{k\rightarrow\infty} \hat N^k(\tau) = N^*(\chi) , \label{eq-4-lem20}
\end{eqnarray}
for some constant $\chi \ge 0$.
Let $\epsilon>0$ be any given (small) number. Then, there exists a
sufficiently large $T$ such that, for sufficiently large $k$,
the following results
 hold for all non-negative integers $j \le k\delta /T$:
\begin{itemize}
\item[(a)]
 (uniform attraction)
 \begin{eqnarray}
  |\bar N^{k,j}(u) - n^*(\bar W^{k,j}(u)) | \le \epsilon ,
  \mbox{~~~~for all } u \in [0, T]; \nonumber
 \end{eqnarray}
\item[(b)]
 (boundedness)
\begin{eqnarray}
  \bar W^{k,j}(u)\le \chi +C+O(\epsilon),
   \mbox{~~~~for all } u \in [0, T]; \nonumber 
\end{eqnarray}
 where $\lim_{\epsilon\to 0} O(\epsilon)=0$;
 i.e.,  $\bar W^{k,j}(u)$ is uniformly bounded;
 and hence, so is $ \bar N^{k,j}(u)$;
\item[(c)]
 (complementarity) if $\bar W^{k,j}(u) >\epsilon$ for
 all $u\in [0,T]$, then
 \begin{eqnarray}
  \bar Y^{k,j}(u) - \bar Y^{k,j}(0)= 0 ,
  \mbox{~~~~for all } u \in [0, T]. \nonumber 
 \end{eqnarray}
\end{itemize}
}
\end{lem}

\medskip

The proof of the lemma is postponed to the appendix.
Below, we first comment on the key points of the lemma, followed
by proving Theorem \ref{thm-opt-diff} in the next subsection.

\medskip\noindent
{\bf Remarks.}
\begin{itemize}
 \item
The main difficulty in proving Theorem \ref{thm-opt-diff} is the
lack of complementarity between $\hat W^k(t)$ and $\hat Y^k(t)$.
Specifically, while (\ref{eq-diff-110}) and (\ref{eq-diff-146}) hold
for $\hat W^k(t)$ and $\hat Y^k(t)$ as counterparts to
(\ref{eq-RBM-10}) and (\ref{eq-RBM-50}) for the limiting processes
$\hat W(t)$ and $\hat Y(t)$, the complementarity in
(\ref{eq-RBM-50a}) does not necessarily hold for $\hat W^k(t)$ and
$\hat Y^k(t)$. Indeed, for an arbitrarily given $k$, it may not hold
when $\hat W^k(t)$ is positive but $\hat N^k(t)$ deviates
significantly from the fixed point state. Part (c) of the above
lemma asserts, however, that complementarity will hold when $k$ is
large enough.
 \item
Recall from Theorem \ref{thm-fluid-ua}, we know that the
fluid-scaled queue-length process will approach the fixed-point
after a certain time of length $T$. This means, without the fluid scaling,
the time to reach the fixed-point is in the order of $kT$, which
translates into a time length in the order of $T/k$ for the
diffusion-scaled process $\hat N^k(t)$. In this sense, $\hat N^k(t)$
approaches the fixed point almost ``instantaneously". This is
confirmed in part (a) of the above lemma, i.e., $\hat N^k(t)$
evolves closely around the fixed-point process (for sufficiently
large $k$). As the fixed point is unique given the workload $w$, the
 process $\hat N^k(t)$ can be characterized
(approximately)
by the workload process $\hat W^k(t)$.
\item
Furthermore, Lemma \ref{lem-y-n} guarantees that,
under the single-bottleneck condition, in the fixed-point
state, there will be no unused capacity at the bottleneck
link, as long as there is some positive workload (no matter
how small) associated with the fixed-point
state. Hence,
 $\hat Y^k(t)$ will not increase, which is the required complementarity.
\end{itemize}

\subsection{Proof of Theorem \ref{thm-opt-diff}}
 \label{sec-diff-proof1}


As before, we focus on a fixed sample path.
Since $\hat Y^k (t)$, for each $k$, is a process that is
nondecreasing and right continuous with left limit (RCLL), we are
guaranteed that for any subsequence of these processes there exists
a further subsequence, denoted $\calk$, that converges to a limit
$\hat Y(t)$, which is also nondecreasing and RCLL; note that $\hat
Y(t)$ is continuous for almost all time $t$ and that {\it at the
moment} this convergence is guaranteed {\it only} for those time $t$
at which $\hat Y(t)$ is continuous. Consequently, we have, $\hat W^k
(t)\to \hat W (t)$, along the same convergent subsequence $\calk$,
with $\hat W(t)$ satisfying the relation in (\ref{eq-RBM-10}): $\hat
W (t)=\hat X (t)+\hat Y (t)$. Furthermore, $\hat W (t)$ is bounded,
following part (b) of Lemma \ref{lem-diff-key}. (We can choose $\tau
= 0$ and any $\delta$ in the lemma; and hence, $\chi=0$.)
 This implies that $\hat Y (t)$ is
also bounded, since $\hat X(t)$ is bounded by $C$ (cf.\
(\ref{eq-4-C-bdd})).

Based on part (a) of Lemma \ref{lem-diff-key},
we have, also  along the subsequence $\calk$,
$\hat N^k (t)\to \hat N(t) := n^*(\hat W(t))$.

Next, we argue that the limit, $\hat Y(t)$ is, in fact, {\it
continuous}; hence, so are $\hat W(t)$ (since $\hat X(t)$ is
continuous) and $\hat N(t)$ (since $n^*(w)$ is continuous according
to Lemma \ref{lem-fixPt-work-continue}). That is, the convergence of
$\hat Y^k(t)$, $\hat W^k(t)$ and $\hat N^k(t)$ to their limits holds
for all time $t$, not just for the time points at which they are
continuous as argued above.

Use contradiction. Suppose $\hat Y(t)$ is discontinuous at
$t_J$, specifically,
\begin{eqnarray}
 && C_J := \hat Y(t_J+) - \hat Y(t_J-) >0 . \label{eq-4-500}
\end{eqnarray}
Pick the time interval $[\tau,\tau+\delta]$ to include $t_J$; where
$\tau$ is a continuous time of $\hat Y(t)$ and close enough to $t_J$
and $\delta$ is small enough, such that (i) the inequality in
(\ref{eq-4-C-bdd}) holds with $C<C_J/2$ and (ii) $|\hat W(\tau) -
\hat W(t_J-)| < C$. (This is possible since $\hat X(t)$ is
continuous and $\hat W(t)$ has left limit at any time $t$.) Based on
(i), we can invoke part (b) of Lemma \ref{lem-diff-key}, letting
$\epsilon\to 0$, and letting $\hat W(\tau)=\chi$ be the limit of
$\hat W^k(\tau)$ along the subsequence $\calk$. This, along with
(ii), yields
$$ \hat W(t) \le \chi +C = \hat W(\tau) + C
< \hat W(t_J-) + 2C.$$
From (\ref{eq-4-500}) and the continuity of $\hat X$, we have
$$\hat W(t_J-)
=\hat X(t_J-)+\hat Y(t_J-)
=\hat X(t_J+)+\hat Y(t_J+)-C_J
=\hat W(t_J+)-C_J.$$
Combining the above,
we have
\begin{eqnarray*}
 \hat W(t) &<& \hat W(t_J-) + 2C   \\
 & = & (\hat W(t_J-) + C_J)  - (C_J - 2C) \\
 & = & \hat W(t_J+) - (C_J - 2C) .
\end{eqnarray*}
That is,
$$\hat W(t_J+)- \hat W(t) > C_J - 2C >0$$
for any $t\in [\tau,\tau+\delta]$, contradicting the
right continuity of $\hat W$ (at $t=t_J$).

We can now make use of part (c) of Lemma
\ref{lem-diff-key} to claim the
 complementarity in (\ref{eq-RBM-50a}), i.e.,
 if $\hat W(t) >0$ for some
time $t>0$, then there exist a small time interval $(\tau,
\tau+\delta)\ni t$
such that  $\hat Y(t)$ does not increase, i.e., $\hat Y(s) - \hat
Y(\tau) = 0$ for all $ s\in (\tau, \tau+\delta)$. The discussion of
the case $t=0$ is similar and hence omitted.

Having proved that the convergence, along the subsequence $\calk$,
to the limit $(\hat Y(t),\hat W(t),\hat N(t))$ holds for all $t$,
and that the limit is continuous and satisfies all the requirements
in (\ref{eq-RBM-10})-(\ref{eq-RBM-fixed}), we can invoke the
uniqueness of the solution to the Skorohod problem (e.g.,
\cite{cynet}, Chapter 6) to conclude that the convergence holds for
the original (full) sequence u.o.c.

Finally, we prove the optimality in part (b) of Theorem \ref{thm-opt-diff}.
Let
\begin{eqnarray}
 \hat Y^G(t) := \liminf_{k\rightarrow\infty} \hat Y^{k,G}(t) ,
 \nd
 \hat W^G(t) := \liminf_{k\rightarrow\infty} \hat W^{k,G}(t);
 \nonumber
\end{eqnarray}
and note that $\hat X^{k,G}(t)\to\hat X(t)$ is independent of
 the resource allocation scheme $G$.
The following is then directly verified:
\begin{eqnarray*}
 && \hat W^G(t) = \hat X(t) + \hat Y^G(t) \ge0 \mbox{~~~~for all } t\ge0 ,\\
 && \hat Y^G(t) \mbox{ is non-decreasing with } \hat Y^G(0) \ge0 ;
\end{eqnarray*}
whereas the complementarity need not hold for $( \hat W^G(t), \hat Y^G(t))$.
Hence, we have
$$ \liminf_{k\to\infty} \hat W^{k,G}(t) \ge \hat W (t),$$
following the well known minimality of the Skorohod problem.
This inequality then implies
$$\liminf_{k\to\infty} \sum_{r\in\calr} C_r(\hat N_r^{k,G}(t))
 \ge \sum_{r\in\calr} C_r(\hat N_r(t)) ,$$
since  $\hat N(t)$ is a fixed-point
(see  (\ref{eq-RBM-fixed})) and hence a minimizer of the cost function.
The above two inequalities  then imply the two inequalities
in (\ref{eq-4-55}).
\hfill$\qed$

\section{Appendix}

\subsection{Proof of Lemma \ref{lem-y-n}}

Let ${\cal R}^* := \{ r\in {\cal R}: \ell^* \in r \}$ denote the set
of routes that involve the bottleneck link $\ell^*$. Partition the
link set $\call$ into ${\cal L} =\{\ell^*\} \cup {\cal L}' \cup
{\cal L}''$, where ${\cal L}' := \{\ell\neq\ell^*, \,\ell\in r
\mbox{ for some } r\in {\cal R}^* \}$ denotes the set of links,
excluding the bottleneck link $\ell^*$, that are used by at least
one route $r\in {\cal R}^*$, and $\call''$ collects all other links.

We prove the lemma by contradiction. Suppose the
lemma does not hold. Then, there is a sequence of states $\tilde n^i
= \{ \tilde n_r^i \ge 0, r\in {\cal R}\}$ ($i=1,2,...$) satisfying
\begin{eqnarray*}
|\tilde n^i - n^*(\tilde w^i)| \rightarrow 0
   \mbox{~~~~as } i\rightarrow \infty ,
\qquad
\tilde  w^i = \sum_{r\ni \ell^*} \nu_r \tilde n_r^i \ge \epsilon,
 \mbox{~~~~for all }
i=1,2, \cdots ,
\end{eqnarray*}
such that
\begin{eqnarray*}
 && \sum_{r\ni\ell^*} \Lambda_r(\tilde n^i) < c_{\ell^*}
 \mbox{~~~~for all } i=1,2, \cdots .
\end{eqnarray*}
We scale the sequence of states by letting $n^i =
\frac{\epsilon}{\tilde w^i} \tilde n^i$, so that $w^i = \sum_{r\ni
\ell^*} \nu_r n_r^i = \epsilon$, for all $i=1,2, \dots$. Then, from
the above
and the homogeneity
property  in (\ref{eq-assume-4}), we have
\begin{eqnarray}
 n^i \rightarrow n^*(\epsilon)
 \mbox{~~~~as } i\rightarrow \infty ;
 \qquad \sum_{r\ni\ell^*} \Lambda_r(n^i) < c_{\ell^*}
 \mbox{~~~~for all } i=1,2, \cdots .
\label{eq-2-y60}
\end{eqnarray}

Rewrite the utility optimization problem
as follows:
\begin{eqnarray*}
 \max && \sum_{r\in {\cal R}^*} U_r(n_r, \Lambda_r)
     + \sum_{r\in {\cal R} \setminus {\cal R}^*} U_r(n_r, \Lambda_r)\\
 s.t. && \sum_{r\in {\cal R}^*: r\ni \ell^*} \Lambda_r
  \phantom{~~~~~~~~+ \sum_{r\in {\cal R}\setminus {\cal R}^*:r\ni\ell} \Lambda_r}
  \le c_{\ell^*} \\
 && \sum_{r\in {\cal R}^*:r\ni\ell} \Lambda_r
  ~~~~+~~ \sum_{r\in {\cal R}\setminus {\cal R}^*:r\ni\ell} \Lambda_r \le c_{\ell}
   \mbox{~~~~for } \ell \in {\cal L}' , \\
 && \phantom{\sum_{r\in {\cal R}^*:r\ni\ell} \Lambda_r ~~~~+~~ }
 \sum_{r\in {\cal R}\setminus {\cal R}^*:r\ni\ell} \Lambda_r \le c_{\ell}
 \mbox{~~~~for } \ell \in {\cal L}'' , \\
 && \Lambda_r \ge 0  \mbox{~~~~for } r \in {\cal R}  .
\end{eqnarray*}
As the allocation $\Lambda(n^i)$ is the optimal solution to the
above, there exists a set of Lagrange multiplies $\{ \theta_\ell^i,
\ell\in {\cal L} \}$ such that, by way of the KKT condition and the
inequality in (\ref{eq-2-y60}), the following holds for all
$i=1,2,\cdots$,
\begin{eqnarray}
 && \theta_\ell^i \ge 0 \mbox{~~~~for } \ell \in {\cal L} , \nonumber\\
 && \theta_{\ell^*}^i = 0 , \nonumber\\
 && \partial_2 U_r(n_r^i, \Lambda_r(n^i))
   \le \sum_{\ell\in r} \theta_\ell^i
   = \sum_{\ell\in r \cap {\cal L}'} \theta_\ell^i
   \mbox{~~~~for } r\in {\cal R}^* ,
   \label{eq-2-y100}\\
 && \partial_2 U_r(n_r^i, \Lambda_r(n^i))
   \le \sum_{\ell\in r} \theta_\ell^i
   = \sum_{\ell\in r \cap {\cal L}'} \theta_\ell^i
        + \sum_{\ell\in r \cap {\cal L}''} \theta_\ell^i
   \mbox{~~~~for } r\in {\cal R}\setminus {\cal }R^* .
   \label{eq-2-y188}
\end{eqnarray}

First, summing up all the inequalities in (\ref{eq-2-y100}) yields
\begin{eqnarray}
 \sum_{r\in {\cal R}^*} \partial_2 U_r(n_r^i, \Lambda_r(n^i))
   &\le& \sum_{r\in {\cal R}^*} \sum_{\ell\in r \cap {\cal L}'} \theta_\ell^i
   \nonumber\\
   & \le & \sum_{\ell\in {\cal L}'} \sum_{r\in {\cal R}^*}  \theta_\ell^i
     \le |{\cal R}|\sum_{\ell\in {\cal L}'}  \theta_\ell^i ,
     \nonumber
\end{eqnarray}
where $ |{\cal R}|$  denotes the total number of routes in the set $\calr$.
Hence, we have
\begin{eqnarray}
 \liminf_{i\rightarrow \infty} \sum_{\ell\in {\cal L}'} \theta_\ell^i
  &\ge& \lim_{i\rightarrow \infty} \frac{1}{|{\cal R}|}
      \sum_{r\in {\cal R}^*} \partial_2 U_r(n_r^i, \Lambda_r(n^i))
\nonumber\\
 & =& \frac{1}{|{\cal R}|} \sum_{r\in {\cal R}^*} \partial_2 U_r(n^*(\epsilon), \rho_r)
  > 0 , \label{eq-2-y120}
\end{eqnarray}
where the equality follows from the
continuity of the allocation $\Lambda(\cdot)$ (Lemma
\ref{lem-lambda-pty}), property (i) in Proposition
\ref{pro:utilcost},
and the continuity of $\partial_2 U_r(\cdot, \cdot)$ (recall
$U_r$ is assumed to be twice-differentiable).

Next,
choose
$\alpha>0$ and $\beta>0$, sufficiently small, such that
\begin{eqnarray}
 && \sum_{r\in {\cal R}^*: r\ni \ell} \rho_r + \alpha
      + |{\cal R}\setminus {\cal R}^*| \beta < c_\ell
  \mbox{~~~~for all } \ell \in {\cal L}' .
  \nonumber
\end{eqnarray}
For each index $i$, denote
${\cal L}'(i) := \{ \ell\in {\cal L}' : \theta_\ell^i>0 \}$.
From the continuity of the allocation, we know $\Lambda_r(n^i)
\rightarrow \Lambda_r(n^*(\epsilon)) = \rho_r$ for all $r\in {\cal
R}^*$. Thus, for sufficiently large $i$, we have
\begin{eqnarray}
 && \sum_{r\in {\cal R}^*: r\ni \ell} \Lambda_r(n^i)
   < \sum_{r\in {\cal R}^*: r\ni \ell} \rho_r + \alpha
   \mbox{~~~~for } \ell \in {\cal L}' .
   \nonumber
\end{eqnarray}
Since $\theta_\ell^i >0$ for any $\ell \in {\cal L}'(i)$, we have,
according to the KKT condition,
\begin{eqnarray}
 && \sum_{r\ni \ell} \Lambda_r(n^i) = c_\ell
   \mbox{~~~~for } \ell \in {\cal L}'(i) .
   \nonumber
\end{eqnarray}
Consequently, for sufficiently large $i$, we have,
 for each $\ell \in {\cal L}'(i)$,
\begin{eqnarray*}
 \sum_{r\in {\cal R}\setminus {\cal R}^* : r\ni \ell} \Lambda_r(n^i)
 & = & c_\ell - \sum_{r\in  {\cal R}^* : r\ni \ell} \Lambda_r(n^i) \\
 & \ge& c_\ell - \left( \sum_{r\in {\cal R}^* : r\ni \ell} \rho_r -\alpha \right)
  > |{\cal R}\setminus {\cal R}^*| \beta .
  \end{eqnarray*}
As the number of terms in the first summation is at most $|{\cal R}
\setminus {\cal R}^*|$, there exists at least one route, denoted as
$r(\ell) \in {\cal R} \setminus {\cal R}^*$ (with $\ell\in {\cal
L}'(i)$), such that
\begin{eqnarray}
 && \Lambda_{r(\ell)} (n^i) > \beta . \label{eq-2-y200}
\end{eqnarray}
According to the KKT condition, any inequality in
(\ref{eq-2-y188}) that corresponds to a route $r(\ell)$, $\ell \in
{\cal L}'(i)$, should hold as equality. Summing up all such equalities
yields
\begin{eqnarray}
 && \sum_{\ell\in {\cal L}'(i)}
    \partial_2 U_{r(\ell)}(n_{r(\ell)}^i, \Lambda_{r(\ell)} (n^i))
  = \sum_{\ell\in {\cal L}'(i)}
   \left( \sum_{j\in r(\ell) \cap {\cal L}'} \theta_j^i
   + \sum_{j\in r(\ell) \cap {\cal L}''} \theta_j^i \right)
   \nonumber \\
 && ~~~~ \ge
   \sum_{\ell\in {\cal L}'(i)} \sum_{j\in r(\ell) \cap {\cal L}'} \theta_j^i
  \ge \sum_{\ell\in {\cal L}'(i)} \theta_\ell^i
  = \sum_{\ell\in {\cal L}'} \theta_\ell^i .
   \label{eq-2-a30}
\end{eqnarray}
where the second inequality is due to the fact that the summation
$\sum_{j\in r(\ell) \cap {\cal L}'} \theta_j^i$ contains the term
$\theta_\ell^i$, and the last equality follows from the definition
of the set ${\cal L}'(i)$.

On the other hand, in view of (\ref{eq-2-y200}), we have
\begin{eqnarray}
 && \sum_{\ell\in {\cal L}'(i)} \partial_2 U_{r(\ell)}(n_{r(\ell)}^i, \Lambda_{r(\ell)} (n^i))
  \le \sum_{\ell\in {\cal L}'(i)} \partial_2 U_{r(\ell)}(n_{r(\ell)}^i, \beta
  ) \nonumber\\
 && \le \sum_{\ell\in {\cal L}'(i)} \sum_{r\in {\cal R}\setminus {\cal R}^*}
        \partial_2 U_r(n_r^i,\beta)
  \le |{\cal L}'| \sum_{r\in {\cal R}\setminus {\cal R}^*} \partial_2 U_r(n_r^i,\beta)
  \nonumber\\
 && \rightarrow 0 \mbox{~~~~as } i\rightarrow \infty . \label{eq-2-y370}
\end{eqnarray}
(Note that, should the set ${\cal L}'(i)$ be empty for some index
$i$, the above
still holds by default.) Putting together
 (\ref{eq-2-a30}) and (\ref{eq-2-y370}),
we have
$\sum_{\ell\in {\cal L}'} \theta_\ell^i \rightarrow 0$,
as $i\rightarrow \infty$,
contradicting the limit in  (\ref{eq-2-y120}).
 \hfill $\qed$

\subsection{Proof of Lemma \ref{lem-diff-key}}

\subsubsection*{Preparations}

Let $b_N$ and $b_W$ be any constants such that the following
inequalities hold,
\begin{eqnarray*}
 w \le b_N |n^*(w)|, \qquad
 |n^*(w)| \le b_W w.
\end{eqnarray*}
In fact, since
\begin{eqnarray*}
\left(\min_{r\ni \ell^*} \nu_r \right) \sum_{r\ni \ell^*} n_r^*
\le w= \sum_{r\ni \ell^*} \nu_r n_r^*
 \le \left(\max_{r\ni \ell^*} \nu_r \right) \sum_{r\ni \ell^*} n_r^* ,
\end{eqnarray*}
the two constants can be chosen as
\begin{eqnarray*}
  b_N= \max_{r\ni \ell^*} \nu_r ,
\qquad b_W=  \left( \min_{r\ni \ell^*} \nu_r \right)^{-1} .
\end{eqnarray*}
Next, define the following constants,
\begin{eqnarray*}
  B_{N,1}=  \eta(b_W\epsilon+\epsilon) + \epsilon  ,
 && B_{W,1}= b_N B_{N,1} + \epsilon ; \\
 B_{W,2}=  \max \{ B_{W,1} , \chi+\epsilon \} + (C+\epsilon),
 && B_{N,2}= b_W B_{W,2} + \epsilon ;
\end{eqnarray*}
where the function $\eta(\cdot)$ appeared in
Lemma \ref{lem-fluid-ua}, and all other quantities are
specified in the statement of the lemma under proof.
The constants defined above
will be used to bound processes $\hat
N^k(t)$ and $\hat W^k(t)$ for
 $t\in [\tau, \tau + \delta]$ and
sufficiently large $k$.

Next, we specify the time length $T$
(stated in the lemma under proof)
as follows:
\begin{eqnarray}
 T & \ge & \max \{
   T_{b_W\epsilon + \epsilon , \epsilon} ,
   T_{\max \{ B_{N,1} , B_{N,2}\}, \epsilon/2},
   T_{\max \{ B_{N,1} , B_{N,2}\},\sigma/2}
 \} , \label{eq-4-T}
\end{eqnarray}
where the terms on the right hand side were used in
the proof of Theorem
\ref{thm-fluid-ua}, and $\sigma=\sigma(\epsilon)$ is
specified in Lemma \ref{lem-y-n}
(with $\epsilon$ given in the lemma under proof).
Note that $T$ is long enough so that
in the fluid network, the fluid state
$\bar N(t)$ will be
close enough (by an error bound of $\epsilon$) to the fixed-point state,
from an initial state $\bar N(0)$ that is bounded by $b_W\epsilon +
\epsilon$;
or close enough by an error bound of $\epsilon/2$ or $\sigma/2$),
from an initial state that is bounded
by $\max \{ B_{N,1} , B_{N,2}\}$.

With the quantities defined or specified above, we state what we
want to prove, in terms of parts (b) and (c) of the lemma, in the
following stronger form (part (a) remains the same): %
For sufficiently large $k$, the following results hold for all
non-negative integers $j \le k\delta /T$:
\begin{itemize}
\item[{\rm (a)~~}]
  $|\bar N^{k,j}(u) - n^*(\bar W^{k,j}(u)) | \le \epsilon ,
  \mbox{~~~~for all } u \in [0, T];$
\item[{\rm (b1)}]
 if $\bar W^{k,j}(u) \le \epsilon (<C)$ for some $u\in [0,T]$,
 then, for all  $u \in [0, T]$,
 \begin{eqnarray}
 \bar W^{k,j}(u) \le  B_{W,1} ,
  \qquad
|\bar N^{k,j}(u)| \le  B_{N,1} ; \label{eq-4-205}
 \end{eqnarray}
\item[{\rm (b2)}]
 if $\bar W^{k,j}(u) >\epsilon$ for all $u\in [0,T]$,
 then,  for all  $u \in [0, T]$,
 \begin{eqnarray}
  \bar W^{k,j}(u) \le  B_{W,2} \left(= \chi+C+O(\epsilon)\right) ,
  \qquad |\bar N^{k,j}(u)| \le  B_{N,2} , \label{eq-4-220}
 \end{eqnarray}
and
 \begin{eqnarray}
   \bar Y^{k,j}(u) - \bar Y^{k,j}(0) = 0 . \label{eq-4-208}
 \end{eqnarray}
\end{itemize}

\subsubsection*{Step 1 of the Proof}

Here we prove the three parts of the lemma, (a, b1, b2), for $j=0$. Note
that by way of the construction, we have
\begin{eqnarray*}
 && (\bar W^{k,0}(0),\bar N^{k,0}(0)) = (\hat W^k(\tau), \hat N^k(\tau)) ,
  \nonumber
\end{eqnarray*}
and hence,
\begin{eqnarray*}
 && (\bar W^{k,0}(0), \bar N^{k,0}(0)) \rightarrow (\chi, n^*(\chi)) ,
  ~~~~ \mbox{ as } k\rightarrow\infty ,
\end{eqnarray*}
following (\ref{eq-4-lem20}). Then, from Theorem
\ref{thm-stoch2fluid} and Corollary \ref{cor-fluid-ua} (with $\bar
W(0)$ and $\bar N(0)$ replaced by $\chi$ and $n^*(\chi)$
respectively),
we have, as $k\rightarrow \infty$,
\begin{eqnarray*}
 (\bar W^{k,0}(u), \bar N^{k,0}(u)) \rightarrow
 (\bar W(u), \bar N(u))
 \qquad {\rm u.o.c.} \mbox{ of } u\in [0,T],
\end{eqnarray*}
 with $(\bar W, \bar N)$ satisfying $( \bar W(u), \bar N(u) ) = (
\chi, n^*(\chi) )$, for  $u\ge0$. (Note that the convergence here is
along the whole sequence of $k$ rather than a subsequence since the
limit is unique.) Since $n^*(w)$ is continuous in $w$ (Lemma
\ref{lem-fixPt-work-continue}), there exists $\sigma'>0$ such that,
for any feasible workload level $w'$ satisfying $|w' - w|\le
\sigma'$, we have,
$$| n^*(w') - n^*(w) | \le \epsilon /2 .$$
 Let $k$ be sufficiently large such that
$$ |\bar W^{k, 0}(u) - \chi | \le \sigma' , \nd
| \bar N^{k, 0}(u) - n^* (\chi)| \le \epsilon /2 ,$$
 for all $u\in [0,T]$.
Then, combining the above, we have,
$$| n^* (\bar W^{k, 0}(u)) - n^*(\chi)|  \le  \frac{\epsilon}{2}
 \mbox{~~~~for all } u \in [0, T] ,$$
and furthermore,
$$| \bar N^{k, 0}(u) - n^* (\bar W^{k, 0}(u))|
\le  | \bar N^{k, 0}(u) - n^*(\chi)|
  + | n^* (\bar W^{k, 0}(u)) - n^*(\chi)|
\le  \epsilon,$$
for all $u \in [0, T]$. That is, (a) holds when $j =0$ for
sufficiently large $k$.

We now verify (b1, b2). First, from
the established result in (a), we know that $\bar W^{k, 0}(u)$ is
arbitrarily close to $\chi$ for all $u\in [0,T]$ when $k$ is
sufficiently large. This fact directly leads to the inequalities in (b1,b2)
for $j=0$.
Next, recall that the number $\sigma=\sigma(\epsilon)$ is chosen
from Lemma \ref{lem-y-n} with $\epsilon$ given in Lemma
\ref{lem-diff-key}.
For sufficiently large $k$, we have, for all $u \in [0,T]$,
\begin{eqnarray*}
 |\bar N^{k,0}(u) - n^*(\chi)| \le  \frac{\sigma}{2} ,
 \qquad
 |n^*(\chi) - n^*(\bar W^{k,0}(u))| \le \frac{\sigma}{2} ,
\end{eqnarray*}
following what's already established above.
Putting the two together yields
$$ |\bar N^{k,0}(u) - n^*(\bar W^{k,0}(u))|
 \le \sigma .$$
 Finally,
we have, for any $u\in [0,T]$,
\begin{eqnarray*}
 \bar Y^{k,0}(u) - \bar Y^{k,0}(0)
 &=&\int_0^u \sum_{r\ni\ell^*} \left(\rho_r-\Lambda_r(\bar N^{k,0}(s))\right)ds
 \nonumber\\ %
 &=& \int_0^u \left( c_{\ell^*} - \sum_{r\ni\ell^*} \Lambda_r(\bar N^{k,0}(s))\right)ds %
    =0 , \label{eq-4-use-y-lem}
\end{eqnarray*}
where the first equality follows from the definitions of the
processes $\bar Y^{k,j}(u)$ and $\hat Y^k(t)$, along with
(\ref{eq-D-Lambda}), and the fact that $\bar N_r^{k,0}(s)>0$ for all
routes $r\ni \ell^*$ and all time $s\in [0,T]$ under the condition in (a);
the second equality follows from the heavy-traffic condition in
(\ref{htcond}); and the last equality from Lemma \ref{lem-y-n}.
Thus, we have shown the complementarity in (\ref{eq-4-208}) of (b2),
for $j=0$.

\subsubsection*{Step 2 of the Proof}

We now extend the above to $j = 1,\dots, k\delta /T$. Suppose, to
the contrary, there exists a subsequence ${\cal K}_1$ of $k$  such
that, for any $k\in {\cal K}_1$, at least one of the properties (a,
b1, b2) does not hold for some integers $j \in [1, k\delta /T]$.
Consequently, {for any $k\in {\cal K}_1$, there exists a smallest
integer, denoted as $j_k$, in the interval $[1, k\delta /T]$ such
that at least one of the properties (a, b1, b2) does not hold.} To
reach a contradiction, it suffices to construct an infinite
subsequence ${\cal K}_2' \subset {\cal K}_1$, such that the desired
properties in (a, b1, b2) hold for $j=j_k$ for sufficiently large $k
\in {\cal K}_2'$. To construct such a sequence, we will first find a
subsequence ${\cal K}_2 \subset {\cal K}_1$ such that the property
(a) holds for $j=j_k$ for sufficiently large $k \in {\cal K}_2$.
Next, we partition ${\cal K}_2$ into two further subsequences,
${\cal K}_2={\cal K}_3 \cup {\cal K}_4$; and show that the
conclusion of (b1) holds for sufficiently large $k \in {\cal K}_3'
\subset {\cal K}_3$,
and that the conclusion of (b2) holds for sufficiently large $k \in
{\cal K}_4$.
Finally, the subsequence ${\cal K}_2' = {\cal K}_3' \cup {\cal K}_4$
is what we need.

From the proof in Step 1, under what's assumed above,
properties (a, b1, b2) hold for $j = 0, ...,  j_k - 1$, $k \in {\cal
K}_1$. Specifically, for $j = j_k - 1$, we have
\begin{eqnarray}
 && |\bar N^{k, j_k -1}(0)| \le \max \{B_{N,1},B_{N,2}\},
    \mbox{~~~~for all } k \in {\cal K}_1 . \nonumber
\end{eqnarray}
 Therefore, the sequence $\{ \bar N^{k, j_k -1}(0), k\in {\cal K}_1
\}$ has a convergent subsequence.
 Then, by
Theorem \ref{thm-stoch2fluid}, there exists a further subsequence
${\cal K}_2 \subset {\cal K}_1$ such that
\begin{eqnarray}
 (\bar W^{k,j_k -1}(u), \bar N^{k,j_k -1}(u))
   \rightarrow (\bar W(u), \bar N(u))
\qquad{\rm u.o.c.} \mbox{ of } u\ge0 ,
 \mbox{ as } k\rightarrow\infty \mbox{ along } {\cal K}_2 ,
 \label{eq-4-conv10}
\end{eqnarray}
with $|\bar N(0)| \le  \max \{B_{N,1},B_{N,2}\}$.
Then, we have
\begin{eqnarray*}
 && | \bar N^{k, j_k -1}(u) - n^* (\bar W^{k, j_k -1}(u))|
  \nonumber \\
 &\le& | \bar N^{k, j_k -1}(u) - \bar N(u)|
       + | \bar N(u) - n^* (\bar W(u))|
  + | n^* (\bar W(u)) - n^* (\bar W^{k, j_k -1}(u))|
 \nonumber \\
 &\to& | \bar N(u) - n^*(\bar W(u))|
\qquad {\rm u.o.c.} \mbox{ of } u\ge0 , \;   \mbox{ as }
k\rightarrow\infty \mbox{ along } {\cal K}_2 .
   \nonumber
\end{eqnarray*}
Moreover, since $T\ge T_{\max \{B_{N,1},B_{N,2}\},
\epsilon/2}$ and taking into account Corollary \ref{cor-fluid-ua},
we have
\begin{eqnarray}
 &&| \bar N(u) - n^*(\bar W(u))|
   < \frac{\epsilon}{2} \mbox{~~~~for all } u\ge T . \nonumber
\end{eqnarray}
Therefore, for sufficiently large $k\in {\cal K}_2$, we have, for $u
\in [0, T]$,
\begin{eqnarray*}
 && | \bar N^{k, j_k}(u) - n^* (\bar W^{k, j_k}(u))|
    = | \bar N^{k, j_k -1}(T+u) - n^* (\bar W^{k, j_k -1}(T+u))|
 < \epsilon  .
\end{eqnarray*}
That is, (a) holds with $j=j_k$ for sufficiently large $k\in {\cal
K}_2$ ($\subset {\cal K}_1$).

Next, we partition  ${\cal K}_2$ into  ${\cal K}_3 \cup {\cal K}_4$
according to the conditions given in (b1,b2), i.e.,
\begin{eqnarray}
 && {\cal K}_3 = \{ k\in {\cal K}_2:
    \bar W^{k,j_k}(u) \le \epsilon \mbox{ for some } u\in [0,T] \} ,
    \nonumber \\
 && {\cal K}_4 = \{ k\in {\cal K}_2:
    \bar W^{k,j_k}(u) > \epsilon \mbox{ for all } u\in [0,T] \} .
    \nonumber
\end{eqnarray}
Note that at least one of the two sequences ${\cal K}_3$ and ${\cal
K}_4$ must be infinite.

Suppose ${\cal K}_3$ is infinite. Then, for each $k\in {\cal K}_3$,
there exists a fixed $u_k \in [0,T]$ satisfying
\begin{eqnarray}
 &&  \bar W^{k, j_k}(u_k) \le \epsilon . \label{eq-4-ii-20}
\end{eqnarray}
Furthermore, we can choose a subsequence ${\cal K}_3' \subset {\cal
K}_3$ such that, for some $u'\in [0,T]$,
\begin{eqnarray}
 &&  u_k \rightarrow  u'
 \mbox{~~~~as } k\rightarrow\infty \mbox{ along } {\cal K}_3' .
 \nonumber
\end{eqnarray}
Note that the convergence in (\ref{eq-4-conv10})
is valid for the subsequence ${\cal K}_3' \subset {\cal K}_2$ too.
 Before proceeding, we refine the bound for the initial state
$\bar N(0)$. First, note that
\begin{eqnarray}
 && \bar W(0) \le \bar W(u')
 = \lim_{k\rightarrow\infty} \bar W^{k,j_k}(u_k) \le \epsilon ,
 \nonumber
\end{eqnarray}
where the first inequality follows from the increasing property of
$\bar W(t)$;
and the second one from (\ref{eq-4-ii-20}).
 Next, we have the following estimation,
\begin{eqnarray*}
 && |\bar N(0) - n^*(\bar W(0))| \nonumber\\
 &\le& |\bar N(0) - \bar N^{k,j_k}(0)|
   + |\bar N^{k,j_k}(0) - n^*(\bar W^{k,j_k}(0))|
 + |n^*(\bar W^{k,j_k}(0)) - n^*(\bar W(0))| .
\end{eqnarray*}
The first and the third terms on the right hand side above
will approach $0$ as $k\rightarrow\infty$; whereas the second term
is bounded by $\epsilon$ for sufficiently large $k\in {\cal K}_3'$,
since
 property (a) holds for $k\in {\cal K}_3'
\subset {\cal K}_2$.
Hence,
\begin{eqnarray}
 && |\bar N(0) - n^*(\bar W(0))| \le \epsilon ; \nonumber
\end{eqnarray}
from which we obtain a refined bound for  $\bar
N(0)$:
\begin{eqnarray}
 && |\bar N(0)| \le |n^*(\bar W(0))| + \epsilon
   \le b_W \epsilon + \epsilon . \nonumber
\end{eqnarray}
As a by-product, we have
\begin{eqnarray*}
 |\bar N(u)|  \le   \eta(|\bar N(0)|)
  \le \eta(b_W \epsilon + \epsilon) ,
\end{eqnarray*}
where the first inequality follows from Lemma \ref{lem-fluid-ua};
and, letting $u\to\infty$,
\begin{eqnarray*}
 |n^*(\bar W(\infty))| = |\bar N(\infty)|
   \le \eta(b_W \epsilon + \epsilon) \le B_{N,1} ,
\end{eqnarray*}
where the equality follows from the fact that $\bar N(\infty)$ is a
fixed point state (cf. Theorem \ref{thm-fluid-ua}).

 Now, for sufficiently large $k\in {\cal K}_3'$, we have, for all
$u\in [0,T]$,
\begin{eqnarray}
 |\bar N^{k,j_k}(u)| & \le & |\bar N(u)| + \epsilon
  \le  \eta(b_W \epsilon + \epsilon) + \epsilon = B_{N,1} ,
  \label{eq-4-ii-70} \\
 \bar W^{k,j_k}(u) & \le & \bar W(u) + \epsilon
  \le \bar W(\infty) + \epsilon  \nonumber \\
 & \le & b_N |n^*(\bar W(\infty))| + \epsilon
 \le b_N B_{N,1} + \epsilon = B_{W,1} , \label{eq-4-ii-80}
\end{eqnarray}
where
 the first inequality in (\ref{eq-4-ii-70}) follows from
(\ref{eq-4-conv10}),
 and so is the first inequality in (\ref{eq-4-ii-80});
while the second inequality in (\ref{eq-4-ii-80}) follows from the increasing
property of the workload $\bar W$.
The two inequalities in (\ref{eq-4-ii-70}) and (\ref{eq-4-ii-80})
together imply that (b1) holds for $j=j_k$ for sufficiently large
$k\in {\cal K}_3'$.

Next, suppose ${\cal K}_4$ is infinite. The convergence in
(\ref{eq-4-conv10}) is
valid for the subsequence ${\cal K}_4 \subset {\cal K}_2$ too.
Recall, $\sigma = \sigma(\epsilon)$ is chosen from Lemma
\ref{lem-y-n}.
 Then, given the time $T\ge T_{\max \{ B_{N,1} , B_{N,2}
\},\sigma/2}$, we have, for all $u\in [0,T]$,
\begin{eqnarray}
 &&  |\bar N(T+u) - n^*(\bar W(T+u))| \le \frac{\sigma}{2} ,
 \nonumber
\end{eqnarray}
due to Corollary \ref{cor-fluid-ua}.
In addition,  as $k\rightarrow \infty$ along the sequence ${\cal
K}_4$ ($\subset {\cal K}_2$), we have the following convergence
(u.o.c.\ of $u\ge0$),
\begin{eqnarray}
  |\bar N^{k,j_k-1}(u) -\bar N(u)| \rightarrow 0 ,
 \nd |n^*(\bar W(u)) - n^*(\bar W^{k,j_k-1}(u))| \rightarrow 0 .
 \nonumber
\end{eqnarray}
Therefore, for sufficiently large $k\in {\cal K}_4$, we have the
following bound for all $u\in [0,T]$,
\begin{eqnarray*}
 && |\bar N^{k,j_k-1}(T+u) - n^*(\bar W^{k,j_k-1}(T+u))| \\
 &&~~ \le |\bar N^{k,j_k-1}(T+u) -\bar N(T+u)|
  + |\bar N(T+u) - n^*(\bar W(T+u))| \\
 &&~~~~~~ + |n^*(\bar W(T+u)) - n^*(\bar W^{k,j_k-1}(T+u))| \\
 &&~~ \le \sigma .
\end{eqnarray*}
From Lemma \ref{lem-y-n} and the fact that $\bar W^{k,j_k}(u)>
\epsilon$ for all $u\in [0,T]$, we have
\begin{eqnarray}
 && \bar Y^{k,j_k-1}(T+u) \mbox{ does not increase in } u\in [0, T],
 \nonumber
\end{eqnarray}
or equivalently,
\begin{eqnarray}
 && \bar Y^{k,j_k}(u) - \bar Y^{k,j_k}(0) =0
  \mbox{~~~~for all } u\in [0, T] , \label{eq-4-ii-120}
\end{eqnarray}
for sufficiently large $k\in {\cal K}_4$.

Using the complementarity property just established, we estimate the
upper bounds for $\bar W^{k,j_k}(u)$ and $\bar N^{k,j_k}(u)$, for
$u\in [0,T]$.
For a given (sufficiently large) $k \in {\cal K}_4$, there are two
mutually exclusive cases: 1) the condition (as well as the
conclusions) in (b2) holds for all $j = 0, ..., j_k$;
 2) the condition in  (b1) holds for some $j = 0 \le j \le j_k-1$.

 In the first case, the process $\bar Y^{k,j}(u)$ does not
increase in $u\in [0,T]$, for $j = 0, ..., j_k-1$. Thus, we have,
for sufficiently large $k\in {\cal K}_4$,
\begin{eqnarray*}
 \bar W^{k, j_k}(u) &=&  \bar W^{k, 0}(0)
 + \sum_{j=0}^{j_k-1} \left( \bar W^{k, j}(T) - \bar W^{k, j}(0) \right)
 + \left( \bar W^{k, j_k}(u) - \bar W^{k, j_k}(0) \right) \\
 & = & \bar W^{k, 0}(0)
  + \sum_{j=0}^{j_k-1} \left( \bar X^{k, j}(T) - \bar X^{k, j}(0) \right)
  + \left( \bar X^{k, j_k}(u) - \bar X^{k, j_k}(0) \right) \\
 && \phantom{\bar W^{k, 0}(0)}
  + \sum_{j=0}^{j_k-1} \left( \bar Y^{k, j}(T) - \bar Y^{k, j}(0) \right)
  + \left( \bar Y^{k, j_k}(u) - \bar Y^{k, j_k}(0) \right) \\
 & = & \bar W^{k, 0}(0)
  + \sum_{j=0}^{j_k-1} \left( \bar X^{k, j}(T) - \bar X^{k, j}(0) \right)
  + \left( \bar X^{k, j_k}(u) - \bar X^{k, j_k}(0) \right) \\
 & = & \hat W^k(\tau)
  + \left( \hat X^{k}(\tau + j_k T/k + u/k) - \hat X^{k}(\tau) \right) \\
 & \le & (\chi + \epsilon) + (C+\epsilon) ,
\end{eqnarray*}
where the inequality follows from
(\ref{eq-4-C-bdd}) and (\ref{eq-4-lem20}).

Under Case 2, let $j_k^0$ be the largest integer
 such that the condition in (b1) holds. Thus, for all
$j=j_k^0 +1 \le j \le j_k$, the condition and results in (b2) hold,
and hence $\bar Y^{k,j}(u)$ does not increase in $u\in [0,T]$. Then,
similar to Case 1, we have, for sufficiently large $k\in {\cal
K}_4$,
\begin{eqnarray*}
 \bar W^{k, j_k}(u) &=&  \bar W^{k, j_k^0}(T)
 + \sum_{j= j_k^0 +1}^{j_k-1} \left( \bar W^{k, j}(T) - \bar W^{k, j}(0) \right)
 + \left( \bar W^{k, j_k}(u) - \bar W^{k, j_k}(0) \right) \\
 & = & \bar W^{k, j_k^0}(T)
   + \left( \hat X^{k}(\tau + j_k T/k + u/k)
   - \hat X^{k}(\tau + j_k^0 T/k + T/k) \right) \\
 & \le & B_{W,1} + (C+\epsilon) .
\end{eqnarray*}
where the inequality is due to the bound (\ref{eq-4-205}) in (b1)
with $j=j_k^0$ and the definition of the constant $C$ in
(\ref{eq-4-C-bdd}).
 Then, synthesizing the bounds in the two cases, we have, for
sufficiently large $k\in {\cal K}_4$ and for all $u\in [0,T]$,
\begin{eqnarray}
 && \bar W^{k,j_k}(u)
  \le \max \{ (\chi + \epsilon) + (C+\epsilon) ,
       B_{W,1} + (C+\epsilon) \}
  = B_{W,2} ; \nonumber
\end{eqnarray}
and furthermore
\begin{eqnarray*}
 && |\bar N^{k, j_k}(u)|
 \le |n^*(\bar W^{k, j_k }(u))| + \epsilon \\ 
 &&~~~~ \le b_W \bar W^{k, j_k}(u) + \epsilon
 \le b_W B_{W,2} + \epsilon   
   = B_{N,2} .
\end{eqnarray*}
The above two bounds, together with the complementarity property in
(\ref{eq-4-ii-120}), imply that (b2) holds with $j=j_k$ for
sufficiently large $k\in {\cal K}_4$.

Finally, let ${\cal K}_2' = {\cal K}_3' \cup {\cal K}_4$ ($\subset
{\cal K}_2 \subset {\cal K}_1$). Then, the properties in (a, b1, b2)
with $j=j_k$ hold for sufficiently large $k\in {\cal K}_2'$
($\subset {\cal K}_1$).
\hfill $\qed$


\end{document}